\documentclass[11pt]{article}
\usepackage[a4paper,margin=1in]{geometry}
\usepackage{iftex}
\ifPDFTeX
  \usepackage[T1]{fontenc}
  \usepackage[utf8]{inputenc}
  \usepackage{tgtermes,tgheros,tgcursor,comment}
\else
\usepackage{fontspec}
\fi
\usepackage{amsmath,amssymb,amsthm,mathtools}
\usepackage{booktabs,array,longtable,graphicx,tabularx}
\usepackage{microtype,enumitem}
\usepackage[hypertexnames=false]{hyperref}
\hypersetup{colorlinks=true,linkcolor=blue!50!black,citecolor=blue!50!black,urlcolor=blue!50!black,
 pdftitle={A New Computational Model for the Spin-Free Sum-of-Squares Energy Bounds of Electronic Hamiltonian},pdfauthor={Yuming Bao, Yancheng Wu, Dongdong Ge, Yinyu Ye}}
\usepackage{xcolor}
\newcommand{\R}{\mathbb{R}}

\newcommand{\Tr}{\operatorname{Tr}}
\newcommand{\svec}{\operatorname{svec}}
\newcommand{\smat}{\operatorname{smat}}
\newcommand{\vecop}{\operatorname{vec}}

\newtheorem{theorem}{Theorem}[section]
\newtheorem{proposition}[theorem]{Proposition}

\title{A New Computational Model for the Spin-Free Sum-of-Squares Energy Bounds of Electronic Hamiltonian}
\author{%
  \begin{minipage}{0.95\textwidth}
    \centering
    Yuming Bao$^{1,\dagger}$ \quad
    Yancheng Wu$^{2,\dagger}$ \quad
    Dongdong Ge$^{2,*}$ \quad
    Yinyu Ye$^{2}$\\[0.6em]
    {\small
      $^{1}$ Shanghai University of Finance and Economics, China\\
      $^{2}$ Shanghai Jiao Tong University, China\\[0.5em]
      $^{\dagger}$ These authors contributed equally to this
      work and share first authorship.\\
      $^{*}$ Corresponding author ddge@sjtu.edu.cn.
    }
  \end{minipage}%
}
\date{}
\begin{document}
\maketitle
\begin{abstract}
Optimizing spin-free sum-of-squares (SOS) energy bounds for electronic
Hamiltonians leads to semidefinite programs whose expanded size grows
rapidly with the orbital dimension. We derive an exact reformulation of
an established spin-free SOS model as a convex quadratic problem over
a single symmetric matrix constrained between zero and the identity.
Explicit recovery maps connect feasible quadratic-program iterates to
certificates in the original SOS formulation, and identify the
Frank--Wolfe gap with the separation between recovered primal and dual
objectives. We implement the quadratic action through GPU-accelerated
Coulomb contractions and solve the reformulated problem using a cached
implementation of fully corrective Frank--Wolfe, referred to here
as the dimension-reduced simplicial method (DRSM) that caches operator products and optimizes convex-combination
weights in a small corrective quadratic program. Numerical experiments
cover 65 FCIDUMP models, molecular geometries, and regular hydrogen
rings, including original-coordinate recovery checks on nine instances
and comparisons with FISTA and a boundary-point SDP solver.
On eight H100 GPUs, a cold-start calculation for C$_{720}$ with 10,080
spatial orbitals reaches a freshly recomputed Frank--Wolfe gap of
$8.82\times10^{-9}$ hartree in 5.02 hours. Full-space hydrogen-ring
calculations reach 10,000 orbitals, with cyclic reduction providing
further computational savings. These results show how exact
reformulation, data structure exploitation and reuse of Coulomb actions make optimization
of the selected large-scale SOS relaxation practical.

\end{abstract}
\section{Introduction}
The representation of an electronic Hamiltonian influences both the
energy bounds that can be extracted from it and the algorithms used to
simulate it. A sum-of-squares (SOS) representation,
\begin{equation}
 \mathcal H-E_{\mathrm{SOS}}I=\sum_\alpha O_\alpha^\dagger O_\alpha,
 \label{eq:proposed-sos}
\end{equation}
provides an energy lower bound through operator positivity. For a
prescribed generating space for the operators $O_\alpha$, maximizing
this bound leads to a semidefinite program (SDP)~\cite{Ye2021}. The computational
challenge is to optimize the chosen SOS representation at the orbital
dimensions required by electronic-structure applications.

Low et al.~\cite{low2025} connected this task to quantum simulation
through spectral amplification and studied electronic-structure SOS
representations using GPU-accelerated SDP calculations with
cuLoRADS~\cite{han2024}. Rubin, Low, and DePrince~\cite{rubin2026}
developed a broader weighted-SOS framework and its connection to
variational two-particle reduced density matrices, including
spin-adapted constructions and symmetry constraints. Their analysis
also gives a smaller equivalent spin-free generator representation.
These developments motivate examining the optimization structure of
the spin-free model: how much of the expanded SDP can be eliminated
while preserving its optimal bound and recovering its certificates?

For the coefficient layout considered here, the expanded formulation
contains $N^4+2N^2$ scalar equalities, where $N$ is the number of spatial
orbitals. We exploit the symmetric orbital-pair range and the range
conditions imposed by positive semidefiniteness to eliminate redundant
coordinates. Convex duality then gives the spectral-box problem
\begin{equation}
 \min_{0\preceq\Gamma\preceq I_N}
 \phi(\Gamma)
 =4\gamma^T H_2\gamma+2\operatorname{Tr}(H_1\Gamma),
 \qquad \gamma=\operatorname{svec}(\Gamma),
 \label{eq:proposed-box}
\end{equation}
where $\operatorname{svec}$ is the isometric symmetric vectorization
and $H_2\succeq0$ is the fixed two-electron coefficient operator.
The optimal total SOS energy is the minimum in
Eq.~\eqref{eq:proposed-box} plus the core-energy constant.
The optimization variable is a single symmetric $N$-by-$N$ matrix.

The reformulation also provides a recovery map at every feasible
iterate. A spectral-box matrix determines a feasible SOS certificate
and a dual witness in the original SDP coordinates. Their objective
separation equals the Frank--Wolfe gap of the reduced QP. This identity
connects the stopping criterion of the numerical solver to an
objective interval for the original SOS problem, and makes the
reconstruction useful before exact optimality is reached.

The reduced variable has quadratic size, but explicit storage of
$H_2$ still grows quartically with $N$. We address this cost by
expressing its action as a Coulomb contraction. For geometry-based
inputs, an operator interface combines GPU4PySCF's direct integral
backend~\cite{li2024gpu4pyscf} with fixed basis transformations and
symmetric-coordinate mappings. PySCF~\cite{sun2020pyscf} supplies the
molecular basis and one-electron integrals used in model preparation.
Smaller instances can instead use explicitly stored coefficient blocks
through the same Hessian-action interface.

Repeated Coulomb builds can remain expensive even when the coefficient
matrix is not stored. We therefore use a cached adaptation of the fully corrective
Frank--Wolfe scheme described in Algorithm~4 of
Jaggi~\cite{jaggi2013}. This implementation, referred to here as
the dimension-reduced simplicial method (DRSM), uses approximate
corrective solves over retained feasible points. A spectral linear minimization oracle
generates feasible points, and their operator products are retained.
The corrective step optimizes convex-combination weights using a small
quadratic program assembled from these products. Inner weight updates
reuse the cache, while a new feasible point requires a new operator
action. The global Frank--Wolfe gap determines termination of the
full spectral-box problem.

The contributions are organized around three results:
\begin{enumerate}[leftmargin=*]
\item An exact spectral-box convex QP reformulation of the selected spin-free
SOS model, with feasible-point recovery and an identity relating the
Frank--Wolfe gap to the recovered SOS objective interval.
\item A scalable implementation of fully corrective Frank--Wolfe for
the reformulated QP, combining explicit or Coulomb-based Hessian
actions with cached corrective optimization so that repeated
weight updates reuse previously computed operator products.
\item Numerical evaluation in both reduced and original coordinates,
including comparisons across optimization methods, operator
representations, and symmetry reductions.
\end{enumerate}

The experiments cover 65 distinct FCIDUMP models, 11 molecular
geometries, C$_{720}$, and 19 regular hydrogen rings.
For nine FCIDUMP instances, we reconstruct the five original SOS Gram
matrices and check the coefficient equations and positive-semidefinite
constraints using the original integral inputs. We compare DRSM with
matrix-free FISTA on FCIDUMP and molecular instances. The hydrogen-ring
experiments compare full-space and cyclic-reduced DRSM with a
boundary-point SDP solver applied to the bordered formulation.
A cold-start C$_{720}$ calculation with 10,080 orbitals reaches a
Frank--Wolfe gap of $8.82\times10^{-9}$ Ha in 5.02 hours on eight H100
GPUs, and the full-space hydrogen calculations reach 10,000 orbitals in seconds by data structure exploitation.

\paragraph{Relation to existing methods.}
SOS hierarchies provide a general framework for polynomial and
noncommutative optimization~\cite{Lasserre2001,parrilo2003,helton2002,pironio2010}.
In electronic structure, the v2RDM literature has developed large SDP
formulations and specialized solution
methods~\cite{mazziotti2004,mazziotti2011,deprince2024}.
We build on the spin-free generators and particle--hole identities
of the electronic-structure SOS
formulations~\cite{low2025,rubin2026}, and develop the
optimization-level elimination and recovery for the specified model.
This structural approach complements low-rank SDP
methods~\cite{burer2003,han2024} and symmetry-based
reductions~\cite{gatermannparrilo2004}.
DRSM follows simplicial decomposition~\cite{hearn1985simplicial}
and corrective Frank--Wolfe
methods~\cite{lacostejulien2015variants}. 
More specifically, its point-generation and convex-hull correction
steps adapt the fully corrective scheme in Algorithm~4 of
Jaggi~\cite{jaggi2013} to the spectral-box QP, with cached operator
products and approximate correction over retained feasible points.
Recent work by Halbey et al.~\cite{halbey2025quadratic} develops
efficient quadratic corrections for Frank--Wolfe algorithms. Dimension-reduced second-order methods such as
DRSOM~\cite{zhang2023} exploit curvature information
through low-dimensional quadratic subproblems, providing
another related perspective on our corrective optimization.
Here, the corrective optimization is organized around cached Coulomb
actions, the spectral-box oracle, and SOS certificate recovery.

The paper first defines the SOS model and establishes the
reformulation and recovery maps. It then develops the operator
implementation and solvers, followed by numerical experiments.
Detailed proofs and complete benchmark tables are collected in
the appendices.
\section{Electronic Hamiltonian and the selected SOS model}
\label{sec:electronic hamiltonian and the original sos sdp}
\label{sec:model}

\subsection{Spin-free notation and coefficient matrices}
\label{sec:notation}

Let $N\geq1$ denote the number of spatial orbitals, and let
$a_{i\sigma}^{\dagger}$ and $a_{i\sigma}$ be fermionic creation and
annihilation operators, with $\sigma\in\{\uparrow,\downarrow\}$.
Their canonical anticommutation relations are
\begin{equation}
 \begin{aligned}
 a_{i\sigma}a_{j\tau}^{\dagger}
 +a_{j\tau}^{\dagger}a_{i\sigma}&=\delta_{ij}\delta_{\sigma\tau}I,\\
 a_{i\sigma}a_{j\tau}+a_{j\tau}a_{i\sigma}&=0,\qquad
 a_{i\sigma}^{\dagger}a_{j\tau}^{\dagger}
 +a_{j\tau}^{\dagger}a_{i\sigma}^{\dagger}=0.
 \end{aligned}
 \label{eq:anticomm}
\end{equation}
Here $I$ is the identity operator on the full fermionic Fock space.
The spin-free one-body operators are
\begin{equation}
 E_{ij}=\sum_{\sigma}a_{i\sigma}^{\dagger}a_{j\sigma},
 \qquad E_{ij}^{\dagger}=E_{ji}.
 \label{eq:epq}
\end{equation}
For real, spin-independent integrals, the state-dependent Hamiltonian
can be written as
\begin{equation}
 H=\sum_{i,j}h^{(1)}_{ij}E_{ij}
  +\frac12\sum_{i,j,k,l}h^{(2)}_{ijkl}E_{ij}E_{kl},
 \label{eq:hamiltonian}
\end{equation}
where
\begin{equation}
 h^{(1)}_{ij}=t_{ij}-\frac12\sum_{k=1}^{N}g_{ikkj},
 \qquad h^{(2)}_{ijkl}=g_{ijkl}.
 \label{eq:h1h2}
\end{equation}
The coefficients $t_{ij}$ and $g_{ijkl}=(ij|kl)$ are the supplied
one-electron and chemists' two-electron integrals in a real orthonormal
spatial-orbital basis. Appendix~\ref{app:derivation of H} derives
\eqref{eq:h1h2} by reordering the conventional electronic Hamiltonian.
We distinguish the state-dependent operator $H$ from the total
Hamiltonian
\begin{equation}
 \mathcal H=E_{\rm ecore}I+H.
 \label{eq:total-hamiltonian}
\end{equation}
The scalar $E_{\rm ecore}$ is the constant supplied with the active-space
Hamiltonian, or the nuclear-repulsion energy for a full-space molecular
calculation. The algebraic reductions
below act on $H$; the corresponding total energies are obtained by
adding $E_{\rm ecore}$.

Write $\mathbb S^n$ and $\mathbb S_+^n$ for real symmetric and PSD
matrices of order $n$, respectively. Their Frobenius inner product
is $\langle V,W\rangle=\Tr(V^TW)$, and $I_n$ is the ordinary
matrix identity. Operators such as $E_{ij}$ are distinguished from
ordinary matrix units throughout.
All ordered orbital pairs use the row-major index
$p(i,j)=(i-1)N+j$. For an ordinary matrix $V$,
\[
 [\vecop(V)]_{p(i,j)}=V_{ij},\qquad
 \operatorname{mat}(\vecop(V))=V.
\]
Let $P\in\R^{N^2\times N^2}$ implement transposition and let
$e=\vecop(I_N)$. Then
\begin{equation}
 P\vecop(V)=\vecop(V^T),\qquad
 P=P^T=P^{-1},\quad Pe=e,\quad e^Te=N.
 \label{eq:pair-permutation}
\end{equation}
For symmetric matrices, let $d=N(N+1)/2$ and define normalized
coordinates in lexicographic order over $i\leq j$ by
\begin{equation}
 [\svec(V)]_{(i,j)}=
 \begin{cases}V_{ii},&i=j,\\\sqrt2\,V_{ij},&i<j,\end{cases}
 \qquad \smat(\svec(V))=V.
 \label{eq:svec-definition}
\end{equation}
They satisfy $\Tr(VW)=\svec(V)^T\svec(W)$ for symmetric $V,W$.
For an orthogonal spectral decomposition
$V=O\operatorname{diag}(\nu_i)O^T$, write
$V_+=O\operatorname{diag}(\max\{\nu_i,0\})O^T$ for its positive part.

The fixed coefficient matrices are
\begin{equation}
 (H_1)_{ij}=h^{(1)}_{ij},\qquad
 (H_{\rm full})_{p(i,j),p(k,l)}=\tfrac12 h^{(2)}_{ijkl}.
 \label{eq:coefficient-matrices}
\end{equation}
The factor $1/2$ in $H_{\rm full}$ is part of the Hamiltonian
normalization. Real-orbital exchange symmetries imply
\begin{equation}
 H_1=H_1^T,\qquad H_{\rm full}=H_{\rm full}^T,\qquad
 PH_{\rm full}=H_{\rm full}P=H_{\rm full}.
 \label{eq:data-assumptions}
\end{equation}
These are the coefficient assumptions for the reduction below.
Positive semidefiniteness of the two-body coefficient matrix will
be identified explicitly as the feasibility condition of the SOS model.

\subsection{The five-block SOS model}
\label{sec:original-to-pm}

We consider the real, unweighted five-block SOS model, with the
operator identity imposed on the full Fock space. The spin-free
generators are
\begin{equation}
 x_{ij}=\sum_\sigma a_{j\sigma}^{\dagger}a_{i\sigma},\qquad
 y_{ij}=\sum_\sigma a_{j\sigma}a_{i\sigma}^{\dagger}.
 \label{eq:xy}
\end{equation}
Stack $x$ and $y$ in the ordered-pair convention above and set
$z=[x;y]$. Their Gram matrix is
\begin{equation}
 G_{\rm SF}=\begin{bmatrix}G&G'\\G''&G'''\end{bmatrix}
 \in\mathbb S_+^{2N^2},\qquad G''=(G')^T.
 \label{eq:gsf}
\end{equation}
Together with $D_\sigma,Q_\sigma\in\mathbb S_+^N$, it defines
\begin{equation}
 \widehat H_{\rm SOS}=z^{\dagger}G_{\rm SF}z
 +\sum_{\sigma,i,j}\left[
 (D_\sigma)_{ij}a_{i\sigma}^{\dagger}a_{j\sigma}
 +(Q_\sigma)_{ij}a_{i\sigma}a_{j\sigma}^{\dagger}\right].
 \label{eq:hsos}
\end{equation}
Factoring each PSD Gram block expresses \eqref{eq:hsos} as a sum
of operator squares. The complete Gram variable is
\begin{equation}
 X=\operatorname{blkdiag}
 (D_\uparrow,D_\downarrow,Q_\uparrow,Q_\downarrow,G_{\rm SF})\succeq0.
 \label{eq:fiveblocks}
\end{equation}
The certificate identity is
\begin{equation}
 \begin{aligned}
 \widehat H_{\rm SOS}&=H+\lambda I
 =H-\varepsilon_{\rm SOS}I=\mathcal H-E_{\rm SOS}I,\\
 \varepsilon_{\rm SOS}&=-\lambda,\qquad
 E_{\rm SOS}=E_{\rm ecore}-\lambda.
 \end{aligned}
 \label{eq:certificate_identity}
\end{equation}
The original optimization problem is therefore
\begin{equation}
 \begin{aligned}
 \lambda_* =\min_{\lambda,D_\sigma,Q_\sigma,G_{\rm SF}}\quad&\lambda\\
 \text{s.t.}\qquad&\widehat H_{\rm SOS}=H+\lambda I,\\
 &D_\sigma,Q_\sigma\succeq0\quad(\sigma=\uparrow,\downarrow),
 \qquad G_{\rm SF}\succeq0.
 \end{aligned}
 \tag{SOS-SDP}\label{eq:original-sos}
\end{equation}
For every feasible certificate, $-\lambda$ bounds the expectation of
$H$ from below in every normalized state. Thus $E_{\rm ecore}-\lambda$
is a lower bound for $\mathcal H$, including its restriction to any
fixed-electron-number sector. The optimized state-dependent and total
SOS bounds are, respectively,
$\varepsilon_{\rm SOS}^{\star}=-\lambda_*$ and
$E_{\rm SOS}^{\star}=E_{\rm ecore}-\lambda_*$.

\section{Exact reformulation and recoverable optimality gaps}
\label{sec:reduction}
The reduction proceeds through three representations: the ordered-pair
Gram matrix, a fixed bordered matrix on symmetric orbital pairs, and a
single matrix in the spectral box. The maps preserve feasible objectives
through the primal reductions; convex duality then identifies the final QP.

\subsection{Contraction of the particle--hole Gram block}
\label{sec:contraction}
The particle--hole identity and the identity-augmented generator algebra
are established ingredients of Ref.~\cite{rubin2026}.
With the generators \(x_{ij},y_{ij}\) defined in
\eqref{eq:xy}, anticommutation gives
\begin{equation}
 y_{ij}=2\delta_{ij}I-x_{ji}.
 \label{eq:yidentity}
\end{equation}
Let \(c=2I\) be an operator, and form the operator columns
\begin{equation}
 z=\begin{bmatrix}x\\y\end{bmatrix},\qquad
 w=\begin{bmatrix}x\\c\end{bmatrix},\qquad
 z=Bw,\qquad
 B=\begin{bmatrix}I_{N^2}&0\\-P&e\end{bmatrix}
 \in\R^{2N^2\times(N^2+1)}.
 \label{eq:Bmap}
\end{equation}
Since the entries of \(z,w\) are operators, their quadratic forms
use an adjoint, even though the Gram coefficients are real:
\begin{equation}
 z^\dagger G_{\rm SF}z=w^\dagger Mw,\qquad
 M=B^TG_{\rm SF}B\succeq0.
 \label{eq:Mcontract}
\end{equation}
Write
\begin{equation}
 M=\begin{bmatrix}A&s\\s^T&q\end{bmatrix}.
 \label{eq:Mpartition}
\end{equation}
For \(G_{\rm SF}=\left[\begin{smallmatrix}G&G'\\G''&G'''\end{smallmatrix}\right]\),
where \(G''=(G')^T\), direct multiplication gives
\begin{equation}
 \begin{aligned}
 A&=G-G'P-PG''+PG'''P,\\
 s&=(G'-PG''')e,\qquad q=e^TG'''e.
 \end{aligned}
 \label{eq:Asq-maps}
\end{equation}
In particular, all three maps are linear in the original Gram matrix.
If \(S=\operatorname{mat}(s)\), expansion of the cross and constant
terms gives
\[
 w^\dagger Mw
 =x^\dagger Ax+2\sum_{ij}(S_{ij}+S_{ji})E_{ij}+4qI.
\]
Thus the coefficient-matching conditions, with the present scaling, are
\begin{equation}
 \begin{aligned}
 A&=H_{\rm full},\\
 D_\sigma-Q_\sigma+2(S+S^T)&=H_1,
       &&\sigma=\uparrow,\downarrow,\\
 \lambda&=4q+\Tr(Q_\uparrow)+\Tr(Q_\downarrow).
 \end{aligned}
 \label{eq:exact-matching}
\end{equation}
Positive semidefiniteness will further imply $S=S^T$ in
Section~\ref{sec:symmetric-pairs}.

\begin{proposition}[Exact coefficient matching]
\label{prop:matching}
Under \eqref{eq:data-assumptions}, fix five real symmetric PSD
Gram blocks with the generator structure of Section~\ref{sec:model}.
The identity
\(\widehat H_{\rm SOS}=H+\lambda I\) on the full Fock space holds
if and only if \eqref{eq:exact-matching} holds.
\end{proposition}
The proof is given in Appendix~\ref{app:matching}.

Consequently the compressed coefficient SDP is
\begin{equation}
 \begin{aligned}
 \min_{s,q,D_\sigma,Q_\sigma}\quad&
       4q+\Tr(Q_\uparrow)+\Tr(Q_\downarrow)\\
 \text{s.t.}\quad&
       \begin{bmatrix}H_{\rm full}&s\\s^T&q\end{bmatrix}\succeq0,\\
 &D_\sigma-Q_\sigma+2\bigl(\operatorname{mat}(s)
             +\operatorname{mat}(s)^T\bigr)=H_1,\\
 &D_\sigma,Q_\sigma\succeq0,\qquad \sigma=\uparrow,\downarrow .
 \end{aligned}
 \tag{PM}\label{eq:ordered-bordered}
\end{equation}
The two-body block is fixed by the coefficient equations.

For recovery to the original coordinates, \(B\) has full column rank.
Fix the left inverse
\begin{equation}
 L=(B^TB)^{-1}B^T,\qquad LB=I_{N^2+1}.
 \label{eq:left-inverse}
\end{equation}
Every \(M\succeq0\) has a PSD lift
\begin{equation}
 G_{\rm SF}=L^TML,\qquad B^TG_{\rm SF}B=M.
 \label{eq:reverseM}
\end{equation}
This recovers the compressed data exactly, although it need not recover
the same original Gram matrix.

\begin{theorem}[Exact contraction to the bordered formulation]
\label{thm:pm-equivalence}
Under \eqref{eq:data-assumptions}, \eqref{eq:original-sos} and
\eqref{eq:ordered-bordered} have the same set of feasible objective
values. The forward map is \eqref{eq:Mcontract}. Conversely, every
feasible bordered matrix lifts by \eqref{eq:reverseM}, with the four
spin blocks unchanged. Both maps preserve the objective and recover
an optimal solution in either direction whenever an optimum exists.
\end{theorem}
The proof is in Appendix~\ref{app:pm}.

Starting from the bordered formulation \eqref{eq:ordered-bordered},
we reduce the orbital-pair coordinates to their symmetric subspace
and combine the two identical spin constraints.

\subsection{Symmetric orbital-pair coordinates}
\label{sec:symmetric-pairs}

Let $d=N(N+1)/2=\dim\mathbb S^N$, and let $e_i$ denote the $i$th
coordinate vector in $\mathbb R^N$. Define
\begin{equation}
 F_{ii}=e_ie_i^T,\qquad
 F_{ij}=\frac{e_ie_j^T+e_je_i^T}{\sqrt2},\quad i<j.
 \label{eq:symmetricbasis}
\end{equation}
Let $U\in\mathbb R^{N^2\times d}$ contain the vectors
$\vecop(F_{ij})$ in the ordering used by $\svec$. Then
\[
 U^TU=I_d,\qquad UU^T=\frac{I_{N^2}+P}{2},\qquad
 U\svec(V)=\vecop(V)\quad(V\in\mathbb S^N).
\]
The coefficient symmetries in \eqref{eq:data-assumptions} give
\begin{equation}
 H_2:=U^TH_{\rm full}U\in\mathbb S^d,\qquad
 H_{\rm full}=UH_2U^T.
 \label{eq:H2}
\end{equation}

Positive semidefiniteness of the bordered matrix in
\eqref{eq:ordered-bordered} forces
\begin{equation}
 s\in\operatorname{range}(H_{\rm full})
 \subseteq\operatorname{range}(U).
 \label{eq:cross-range}
\end{equation}
Thus $s=Ut$ for the unique vector $t=U^Ts\in\mathbb R^d$.
Writing $R=\operatorname{diag}(U,1)$, we obtain
\begin{equation}
 M=R\begin{bmatrix}H_2&t\\t^T&q\end{bmatrix}R^T,
 \qquad R^TR=I_{d+1},
 \label{eq:MwithH2}
\end{equation}
and therefore
\begin{equation}
 M\succeq0
 \quad\Longleftrightarrow\quad
 \begin{bmatrix}H_2&t\\t^T&q\end{bmatrix}\succeq0.
 \label{eq:bordered}
\end{equation}
Moreover, $\operatorname{mat}(s)=\smat(t)$ is symmetric, so the
one-body contribution in \eqref{eq:ordered-bordered} becomes
$4\smat(t)$.

\subsection{Spin averaging and the compact primal problem}
\label{sec:p0}

Define the averaged matrices
\[
 D=\frac{D_\uparrow+D_\downarrow}{2},\qquad
 Q=\frac{Q_\uparrow+Q_\downarrow}{2}.
\]
They are PSD and satisfy the same one-body equation as each spin
channel. Since $2\Tr(Q)=\Tr(Q_\uparrow)+\Tr(Q_\downarrow)$,
averaging preserves the objective. Conversely, duplicating $D,Q$
in the two spin channels recovers the original spin constraints.
This gives
\begin{equation}
 \boxed{\begin{aligned}
 \lambda_* =\min_{t,q,D,Q}\quad&4q+2\Tr(Q)\\
 \text{s.t.}\quad&
 \begin{bmatrix}H_2&t\\t^T&q\end{bmatrix}\succeq0,\\
 &D-Q=H_1-4\smat(t),\qquad D,Q\succeq0,
 \end{aligned}}
 \tag{P0}\label{eq:p0}
\end{equation}
where $t\in\mathbb R^d$, $q\in\mathbb R$, and $D,Q\in\mathbb S^N$.

\begin{theorem}[Exact primal equivalence]
\label{thm:exact-reduction}
Under \eqref{eq:data-assumptions}, problems
\eqref{eq:original-sos}, \eqref{eq:ordered-bordered}, and
\eqref{eq:p0} are equivalent: they have the same set of
feasible objective values.
\end{theorem}

The proof is given in Appendix~\ref{app:reduction}.
The reduction replaces the five PSD blocks of orders
$N,N,N,N,2N^2$ with three blocks of orders $N,N,d+1$.

\subsection{The spectral-box quadratic program}
\label{sec:spectral}

The fixed-block formulation admits a quadratic dual over an orbital-space spectral box.

\paragraph{Derivation.}

Assume $H_2\succeq0$. Let $r=\operatorname{rank}(H_2)$ and choose a full-column-rank
factor $L_2\in\mathbb R^{d\times r}$ such that
\begin{equation}
H_2=L_2L_2^T.
\label{eq:spectral-factorization}
\end{equation}

The bordered PSD constraint in \eqref{eq:p0} requires
$t\in\operatorname{range}(H_2)=\operatorname{range}(L_2)$.
Thus $t=L_2y$ for a unique $y\in\mathbb R^r$, and
\[
\begin{bmatrix}
H_2&t\\
t^T&q
\end{bmatrix}
=
\begin{bmatrix}
L_2&0\\
0&1
\end{bmatrix}
\begin{bmatrix}
I_r&y\\
y^T&q
\end{bmatrix}
\begin{bmatrix}
L_2^T&0\\
0&1
\end{bmatrix}.
\]
When $r=0$, $L_2$ has no columns, $y$ is the unique vector in $\mathbb R^0$,
and $L_2y=0$, $\|y\|^2=0$. Thus the bordered constraint reduces to $t=0$, $q\ge0$.

Since the outer factor has full column rank,
positive semidefiniteness is equivalent to that of the
middle matrix. Applying the Schur complement to its
identity block gives
\begin{equation}
\begin{bmatrix}
H_2&t\\
t^T&q
\end{bmatrix}\succeq0
\quad\Longleftrightarrow\quad
\exists\,y\in\mathbb R^r:
\quad t=L_2y,\qquad q\geq\|y\|_2^2.
\label{eq:generalized-schur}
\end{equation}

For fixed $y,D,Q$, the objective is minimized over $q$
by choosing $q=\|y\|_2^2$.
Consequently, \eqref{eq:p0} has the same optimal value as
\begin{equation}
\begin{aligned}
\lambda_*=\min_{y,D,Q}\quad&
4\|y\|_2^2+2\Tr(Q)\\
\text{s.t.}\quad&
D-Q=H_1-4\smat(L_2y),\\
&D,Q\succeq0,
\end{aligned}
\tag{P}
\label{eq:spectral-factor-primal}
\end{equation}
where $y\in\mathbb R^r$ and $D,Q\in\mathbb S^N$.

Introduce a symmetric Lagrange multiplier $2\Gamma$
for the equality constraint in
\eqref{eq:spectral-factor-primal}, and write
$\gamma=\svec(\Gamma)\in\mathbb R^d$.
The Lagrangian is
\begin{equation}
\begin{aligned}
\mathcal L(y,D,Q;\Gamma)
={}&4\|y\|_2^2+2\Tr(Q)\\
&+2\Tr\!\left[
\Gamma\bigl(D-Q-H_1+4\smat(L_2y)\bigr)
\right]\\
={}&4\|y\|_2^2
+8\gamma^TL_2y
-2\Tr(H_1\Gamma)\\
&+2\Tr(\Gamma D)
+2\Tr\bigl((I_N-\Gamma)Q\bigr).
\end{aligned}
\label{eq:spectral-lagrangian}
\end{equation}
Here we used the Frobenius-isometric property of $\svec$:
\[
\Tr\bigl(\Gamma\smat(L_2y)\bigr)
=\gamma^TL_2y.
\]

The infimum over $D,Q\succeq0$ is finite if and only if
\begin{equation}
\Gamma\succeq0,
\qquad I_N-\Gamma\succeq0.
\label{eq:spectral-dual-feasibility}
\end{equation}
Indeed, a negative eigenvalue of either coefficient matrix
allows the corresponding trace term to tend to $-\infty$.
When \eqref{eq:spectral-dual-feasibility} holds, both trace
terms have infimum zero.

The remaining minimization over $y$ can be evaluated by
completing the square:
\begin{equation}
\begin{aligned}
4\|y\|_2^2+8\gamma^TL_2y
&=4\|y+L_2^T\gamma\|_2^2
  -4\|L_2^T\gamma\|_2^2,\\
y_\Gamma&=-L_2^T\gamma.
\end{aligned}
\label{eq:spectral-y-minimization}
\end{equation}
Using $L_2L_2^T=H_2$, the dual function is therefore
\[
\inf_{\substack{y\in\mathbb R^r\\D,Q\succeq0}}
\mathcal L(y,D,Q;\Gamma)
=
\begin{cases}
-4\gamma^TH_2\gamma
-2\Tr(H_1\Gamma),
&0\preceq\Gamma\preceq I_N,\\
-\infty,
&\text{otherwise}.
\end{cases}
\]
Hence the Lagrange dual of \eqref{eq:spectral-factor-primal}
is
\[
\max_{0\preceq\Gamma\preceq I_N}
\left\{
-4\gamma^TH_2\gamma
-2\Tr(H_1\Gamma)
\right\}.
\]
Define
\[
\phi(\Gamma)
=4\gamma^TH_2\gamma+2\Tr(H_1\Gamma).
\]
Changing the sign of the dual objective yields the convex
spectral-box quadratic program
\begin{equation}
\boxed{
\begin{aligned}
\min_{\Gamma\in\mathbb S^N}\quad&\phi(\Gamma)\\
\text{s.t.}\quad&0\preceq\Gamma\preceq I_N.
\end{aligned}
}
\tag{QP-SDP}
\label{eq:phi}
\end{equation}

\begin{theorem}[Exact spectral-box reformulation]
\label{thm:box-duality}
If $H_1\in\mathbb S^N$ and $H_2\in\mathbb S_+^d$, then
\eqref{eq:p0} and \eqref{eq:phi} attain their optima,
and their optimal values satisfy
\begin{equation}
\lambda_*
=-\min_{0\preceq\Gamma\preceq I_N}\phi(\Gamma).
\label{eq:boxdual}
\end{equation}
\end{theorem}

The proof is given in Appendix~\ref{app:duality}.
It establishes strong duality through a strictly feasible
point of the square-root formulation
\eqref{eq:spectral-factor-primal}.
Optimal recovery follows from
Proposition~\ref{prop:recovery}.

Consequently, the optimized SOS energies are
\begin{equation}
\begin{aligned}
\varepsilon_{\rm SOS}^{\star}
&=-\lambda_*
=\min_{0\preceq\Gamma\preceq I_N}\phi(\Gamma),\\
E_{\rm SOS}^{\star}
&=E_{\rm ecore}-\lambda_*\\
&=E_{\rm ecore}
+\min_{0\preceq\Gamma\preceq I_N}\phi(\Gamma).
\end{aligned}
\label{eq:total-sos-optimum}
\end{equation}

The factor $L_2$ is used to establish equivalence and strong duality. The numerical algorithms below apply $H_2$ directly as an operator; they do not form $L_2$ or a pseudoinverse.

\subsection{Feasible recovery and a global stopping criterion}
\label{sec:recovery-gap}
The recovery map applies before convergence. For a feasible
$0\preceq\Gamma\preceq I_N$, let $\gamma=\svec(\Gamma)$ and set
\begin{equation}
\begin{aligned}
 t_\Gamma&=-H_2\gamma,&q_\Gamma&=\gamma^TH_2\gamma,\\
 K_\Gamma&=H_1+4\smat(H_2\gamma),&
 D_\Gamma&=(K_\Gamma)_+,\quad Q_\Gamma=(-K_\Gamma)_+.
\end{aligned}
\label{eq:feasible-recovery}
\end{equation}
Define the primal and dual objective values
\begin{equation}
 \lambda_P(\Gamma)=4q_\Gamma+2\Tr(Q_\Gamma),\qquad
 \lambda_D(\Gamma)=-\phi(\Gamma).
 \label{eq:recovered-objectives}
\end{equation}
\begin{proposition}[Recovery at a feasible iterate]
\label{prop:anytime-recovery}
Under the assumptions of Theorem~\ref{thm:box-duality},
\eqref{eq:feasible-recovery} is feasible for \eqref{eq:p0} and lifts to
an original SOS certificate. If $G=\nabla\phi(\Gamma)=2K_\Gamma$ and
$\widehat\Gamma$ is the orthogonal projector onto the negative eigenspace
of $G$, then
\begin{equation}
\begin{aligned}
 \lambda_P(\Gamma)-\lambda_D(\Gamma)
 &=\Tr\bigl(G(\Gamma-\widehat\Gamma)\bigr)
 =:g_{\rm FW}(\Gamma),\\
 E_{\rm ecore}+\phi(\Gamma)-g_{\rm FW}(\Gamma)
 &\leq E_{\rm SOS}^{\star}
 \leq E_{\rm ecore}+\phi(\Gamma).
\end{aligned}
\label{eq:energy-interval}
\end{equation}
\end{proposition}
Appendix~\ref{app:proof-anytime-recovery} proves feasibility of the
recovered certificate and derives the energy interval
\eqref{eq:energy-interval} from the primal--dual objective gap.

Thus the same global gap measures optimization error and the separation
between recovered SOS objectives. Its upper endpoint bounds the optimum
of the selected relaxation; the ground-state lower bound is the lower
endpoint. The spectral-box variable does not need a fixed-particle-number
interpretation.
\paragraph{Recovery at optimality.}

An optimal solution of the reduced QP determines an optimal
solution of the original SOS--SDP.
The reconstruction first produces an optimal solution of
\eqref{eq:p0}, which is then lifted to the original SOS
formulation using the objective-preserving construction
in the proof of Theorem~\ref{thm:exact-reduction}.

Let $\Gamma^\star$ be an optimal solution of \eqref{eq:phi},
and set $\gamma^\star=\svec(\Gamma^\star)$.
Define
\begin{equation}
\begin{aligned}
t^\star
&=-H_2\gamma^\star,
&
q^\star
&=(\gamma^\star)^TH_2\gamma^\star,
\\
K^\star
&=H_1-4\smat(t^\star)
 =H_1+4\smat(H_2\gamma^\star),
\\
D^\star
&=(K^\star)_+,
&
Q^\star
&=(-K^\star)_+.
\end{aligned}
\label{eq:primal-recovery}
\end{equation}
Here, $A_+$ denotes the positive spectral part of a
symmetric matrix $A$.
Thus, the recovery of $(t^\star,q^\star,D^\star,Q^\star)$
requires only matrix operations and an eigendecomposition
of $K^\star$.

\begin{proposition}[Optimal primal recovery]
\label{prop:recovery}
Under the assumptions of
Theorem~\ref{thm:box-duality},
let $\Gamma^\star$ be an optimal solution of
\eqref{eq:phi}.
Then the variables
$(t^\star,q^\star,D^\star,Q^\star)$ defined in
\eqref{eq:primal-recovery} form an optimal solution of
\eqref{eq:p0}, with objective value
\begin{equation}
4q^\star+2\Tr(Q^\star)
=
-\phi(\Gamma^\star)
=
\lambda_*.
\label{eq:recovery-objective}
\end{equation}
Their objective-preserving lift from the proof of
Theorem~\ref{thm:exact-reduction} is an optimal
certificate for the original SOS--SDP.
Consequently,
\[
\varepsilon_{\rm SOS}^{\star}
=
\phi(\Gamma^\star),
\qquad
E_{\rm SOS}^{\star}
=
E_{\rm ecore}+\phi(\Gamma^\star).
\]
\end{proposition}

\section{Operator evaluation and cached solution of the spectral-box Convex QP}
\label{sec:operator-solver}

The spectral-box reformulation makes it possible to solve larger
instances with a standard first-order method. Evaluating the objective
and its gradient requires only matrix--vector products with $H_2$.
When $H_2$ is too large to store, these products can be computed by
matrix-free Coulomb contractions. Combined with projection onto the
spectral box, this gives a projected FISTA solver that already handles
larger instances than direct SDP optimization in our experiments.

On larger instances, however, repeatedly applying the quadratic
operator can account for much of the solve time. This cost motivates
the dimension-reduced simplicial method (DRSM). DRSM stores the
operator product for each new feasible point and reuses the stored
products when optimizing over convex combinations of these points.
The resulting weight updates require no additional operator products.

\subsection{Matrix-free evaluation of the quadratic operator}
\label{sec:operator-mf}

The reduced quadratic operator can be evaluated through a Coulomb
contraction. In a real orthonormal orbital basis, define
\begin{equation}
 [\mathcal J(\Gamma)]_{ij}
 =\sum_{k,l=1}^{N}g_{ijkl}\Gamma_{kl}.
 \label{eq:mf-orthogonal-J}
\end{equation}
The coefficient convention
$(H_{\rm full})_{p(i,j),p(k,l)}=g_{ijkl}/2$
and the isometric symmetric-pair embedding give
\begin{equation}
 \boxed{
 H_2\gamma=\tfrac12\svec(\mathcal J(\Gamma)),
 \qquad \gamma=\svec(\Gamma).
 }
 \label{eq:mf-action}
\end{equation}
Indeed, $U\gamma=\vecop(\Gamma)$, so
$H_{\rm full}U\gamma=\vecop(\mathcal J(\Gamma))/2$;
multiplication by $U^T$ yields the identity.
Consequently,
\begin{equation}
 \phi(\Gamma)
 =2\Tr(\Gamma\mathcal J(\Gamma))+2\Tr(H_1\Gamma),
 \qquad
 G:=\nabla\phi(\Gamma)=2H_1+4\mathcal J(\Gamma).
 \label{eq:mf-value-gradient}
\end{equation}
A single Coulomb contraction therefore provides both the objective
and the gradient.

\paragraph{Atomic-orbital representation.}
For geometry-based inputs, let $\{\chi_\mu\}$ denote the
atomic-orbital (AO) basis and $S_{\rm AO}$ its overlap matrix.
PySCF~\cite{sun2020pyscf} constructs the molecular basis and
evaluates the AO overlap, kinetic-energy, and nuclear-attraction
matrices. We form a fixed transformation $C_{\rm AO}$ satisfying
$C_{\rm AO}^T S_{\rm AO} C_{\rm AO}=I_N$.
The reported geometry-based calculations use full symmetric
L\"owdin orthogonalization,
$C_{\rm AO}=S_{\rm AO}^{-1/2}$, with positive-definite
$S_{\rm AO}$.

For a trial matrix $\Gamma$, form
\begin{equation}
 D_{\rm AO}=C_{\rm AO}\Gamma C_{\rm AO}^T,
 \qquad
 [\mathcal J_{\rm AO}(D_{\rm AO})]_{\mu\nu}
 =\sum_{\rho\sigma}
 (\mu\nu|\rho\sigma)(D_{\rm AO})_{\rho\sigma}.
 \label{eq:mf-ao-contraction}
\end{equation}
Transforming the contracted result back gives
\begin{equation}
 \boxed{
 H_2\gamma=\tfrac12\svec\!\left[
 C_{\rm AO}^T\mathcal J_{\rm AO}(D_{\rm AO})C_{\rm AO}
 \right].
 }
 \label{eq:mf-ao-action}
\end{equation}

\paragraph{Direct integral backend.}
We evaluate $\mathcal J_{\rm AO}(D_{\rm AO})$ using
GPU4PySCF's conventional four-center integral
backend~\cite{li2024gpu4pyscf}.
The backend evaluates screened shell-quartet contributions
and contracts them with $D_{\rm AO}$ to accumulate the
Coulomb matrix. This direct evaluation avoids storing the
complete four-index integral tensor and the pair-space
matrix $H_2$.
The reported direct-integral runs use the Coulomb component
returned by the JK backend. Integral screening thresholds
and arithmetic precision are specified in
Section~\ref{sec:protocol}.

\paragraph{Hessian-action interface.}
Our implementation combines the Coulomb build with symmetric
unpacking, basis transformations, and the coefficient
normalization of the reduced QP.
For a packed input vector $x$, the solver interface returns
\[
 \mathcal A x
 =
 4\,\svec\!\left[
 C_{\rm AO}^T
 \mathcal J_{\rm AO}\!\left(
 C_{\rm AO}\smat(x)C_{\rm AO}^T
 \right)
 C_{\rm AO}
 \right],
 \qquad \mathcal A=8H_2.
\]
Thus the interface applies the Hessian of the quadratic
objective directly.

At initialization, we combine $C_{\rm AO}$ with the
transformation to the integral engine's internal basis
representation. Each subsequent action forms the trial
matrix in that representation and calls the corresponding
low-level integral build. This avoids repeating the separate
AO representation conversions at every call.
Matrix workspaces are reused, and GPU arrays are shared
between PyTorch and CuPy through DLPack.

For the eight-GPU direct-integral calculations,
GPU4PySCF distributes the integral tasks and reduces the
partial Coulomb matrices. The surrounding coordinate
transformations are performed on GPU~0.
Our solver uses the resulting Hessian actions for the
spectral oracle and corrective optimization described below;
the primary direction cache and spectral calculations
also reside on GPU~0.

When the coefficient matrix fits in device memory, the same
solver interface instead applies explicitly stored $H_2$
blocks distributed across the GPUs. The operator interface
therefore supports both direct-integral matrix-free
evaluation and explicit matrix--vector multiplication.

\subsection{Projected FISTA}
\label{sec:fista}

First we use projected FISTA with backtracking~\cite{beck2009fista} to
minimize $\phi$ over the spectral box $0\preceq\Gamma\preceq I_N$.
Let $\Gamma_k$ denote the current feasible iterate, $Y_k$ the
extrapolated point, and $a_k$ the FISTA momentum parameter.
We initialize $\Gamma_0=I_N/2$, $Y_0=\Gamma_0$, and $a_0=1$.

At iteration $k$, backtracking selects a curvature parameter $L_k>0$
and gives the projected-gradient candidate
\[
  X_k=\Pi_{[0,I]}\!\left(Y_k-L_k^{-1}\nabla\phi(Y_k)\right).
\]
For a symmetric matrix $Z=O\operatorname{diag}(z_i)O^T$, the
Frobenius projection onto the spectral box is
\begin{equation}
  \Pi_{[0,I]}(Z)
  =O\operatorname{diag}\!\left(\min\{1,\max\{0,z_i\}\}\right)O^T.
  \label{eq:box-projection}
\end{equation}
Standard FISTA accepts $\Gamma_{k+1}=X_k$ and extrapolates from the
two most recent feasible iterates according to
\[
  a_{k+1}=\frac{1+\sqrt{1+4a_k^2}}{2},
  \qquad
  Y_{k+1}
  =\Gamma_{k+1}
   +\frac{a_k-1}{a_{k+1}}\bigl(\Gamma_{k+1}-\Gamma_k\bigr).
\]

Additionally, our baseline adds two restart rules to this iteration.
\begin{enumerate}
  \item \textit{Objective safeguard before acceptance.}
  If $\phi(X_k)$ exceeds $\phi(\Gamma_k)$ by more than the
  floating-point tolerance, we discard $X_k$, set $Y_k=\Gamma_k$
  and $a_k=1$, and repeat the projected-gradient step with
  backtracking from $\Gamma_k$. We accept the recomputed candidate
  as $\Gamma_{k+1}$. Otherwise, we accept the original $X_k$.

  \item \textit{Directional restart after acceptance.}
  Once $\Gamma_{k+1}$ has been accepted, we apply the adaptive
  restart criterion of O'Donoghue and
  Cand\`es~\cite{odonoghue2015adaptive}. If
  \[
    \bigl\langle Y_k-\Gamma_{k+1},\,
                 \Gamma_{k+1}-\Gamma_k\bigr\rangle_F>0,
  \]
  we retain $\Gamma_{k+1}$ and set $a_{k+1}=1$ and
  $Y_{k+1}=\Gamma_{k+1}$. The next iteration thus starts from the
  accepted point without extrapolation. Otherwise, we use the
  standard momentum update above.
\end{enumerate}

In the direct-integral implementation, evaluating the gradient at
$Y_k$ and testing new backtracking candidates can require repeated
Coulomb contractions. The method in the next subsection stores the
operator products of feasible search points and reuses them when
updating their convex combination weights.

\subsection{Dimension-reduced simplicial method (DRSM)}
\label{sec:cached-cg}
Furthermore, we write the QP in packed coordinates as
\begin{equation}
 f(\gamma):=\phi(\smat(\gamma))
 =c^T\gamma+\tfrac12\gamma^T \mathcal A\gamma,
 \qquad \mathcal A=8H_2,\quad c=2\svec(H_1),
 \label{eq:cg-packed-objective}
\end{equation}
where $\mathcal A$ denotes a linear operator and need not be assembled.
Repeated products with $\mathcal A$ can dominate a matrix-free solve. 
We therefore develop the \emph{dimension-reduced simplicial method} (DRSM).
DRSM combines a spectral linear minimization oracle with cached operator
products and approximately fully corrective optimization over convex
combination weights.

It follows the simplicial-decomposition principle of
point generation and restricted master optimization
\cite{hearn1985simplicial}, and the corrective Frank--Wolfe framework
\cite{jaggi2013,lacostejulien2015variants}. This optimization of combinations of stored points is also related to Anderson acceleration~\cite{Walker2011}.
Here, dimension reduction refers to the repeated corrective optimization
in the weight space rather than the full packed matrix space. The spectral
oracle and global gap still operate on the original spectral-box QP.

\paragraph{Global generation of feasible points.}
At the current iterate, compute $G=\nabla\phi(\Gamma)$. If
$G=O\operatorname{diag}(\mu_i)O^T$, the linear minimization problem has
an explicit spectral solution:
\begin{equation}
 \widehat\Gamma\in\arg\min_{0\preceq Z\preceq I_N}\Tr(GZ),
 \qquad
 \widehat\Gamma=O\operatorname{diag}(\mathbf1_{\{\mu_i<0\}})O^T.
 \label{eq:cg-spectral-oracle}
\end{equation}
In the eigenbasis of $G$, the objective depends only on diagonal entries
of $Z$, each in $[0,1]$; the displayed choice attains the sum of the
negative eigenvalues. This oracle searches the entire spectral box,
not only the span or convex hull of cached points.

\paragraph{Correction over the current convex hull.}
Maintain feasible matrices $\Gamma_i$, their packed coordinates
$v_i=\svec(\Gamma_i)$, and the products $w_i=\mathcal A v_i$. With
$\mathsf V=[v_1,\ldots,v_m]$ and $\mathsf W=[w_1,\ldots,w_m]$, set
\begin{equation}
 \gamma=\mathsf V\alpha,\qquad
 \alpha\in\Delta_m:=\{a\in\mathbb R^m:a\ge0,\ \mathbf1^Ta=1\}.
 \label{eq:cg-convex-mixture}
\end{equation}
Convexity of the spectral box preserves feasibility. Linearity gives
$\mathcal A\gamma=\mathsf W\alpha$ and hence the restricted problem
\begin{equation}
 \min_{\alpha\in\Delta_m}\ell^T\alpha+\tfrac12\alpha^T\mathsf K\alpha,
 \qquad \ell=\mathsf V^Tc,\quad \mathsf K=\mathsf V^T\mathsf W.
 \label{eq:cg-restricted-qp}
\end{equation}
The term \emph{simplicial} refers to the simplex $\Delta_m$ of
convex-combination weights. Its image under $\mathsf V$ is the
convex hull of the cached points. The corrective problem has
$m$ weights, with at most $m-1$ independent degrees of freedom,
compared with $d=N(N+1)/2$ packed matrix coordinates.
It minimizes the original quadratic objective over this convex hull
using the cached products $\mathsf W$, so weight updates require
no new integral contractions.

When a new point $\widehat\gamma=\svec(\widehat\Gamma)$
is added, only its product $\mathcal A\widehat\gamma$ is required; the existing
products are retained. In exact arithmetic $\mathsf K$ is symmetric PSD. 
The implementation forms a symmetric $\mathsf K$ by averaging reciprocal cross-products.

We approximately solve the dimension-reduced simplex QP using Newton steps on a
candidate support, with step lengths limited by nonnegativity and a
pairwise mass-transfer fallback. The pairwise direction transfers weight
from an active index with large restricted gradient to one with small
restricted gradient. Zero-weight points are removed without changing
the represented iterate. The previous combination remains feasible after
adding a point, so an exact corrective solve cannot increase the objective.
The restricted QP is solved approximately at each outer iteration,
while termination is determined by the full-problem Frank--Wolfe gap.

\paragraph{Outer iteration and cost.}
Starting from $\Gamma_0=I_N/2$, compute and store its operator product.
Each outer iteration of DRSM then performs:
\begin{enumerate}
 \item Form $\gamma=\mathsf V\alpha$ and
       $G=\smat(c+\mathsf W\alpha)$.
 \item Evaluate the global spectral oracle and the original-QP gap.
       If the gap is sufficiently small, perform fresh final checks.
 \item Otherwise append $(\widehat\gamma,\mathcal A\widehat\gamma)$ to the cache,
       update $\mathsf K$ and $\ell$, and approximately solve the restricted QP.
 \item Remove zero-weight points, record diagnostics, and repeat. Checkpoints follow the schedule of each experiment.
\end{enumerate}
The two FP64 caches occupy $16dm$ bytes, and the corrective QP has
$m$ variables. Each outer iteration also requires a spectral
decomposition of the $N$-by-$N$ gradient matrix.
The computational cost is therefore governed by the Coulomb actions,
the spectral decompositions, and the number of retained points.

\subsection{Connection to the recovered SOS interval}
\label{sec:cg-verification}
The stopping quantity is precisely the gap already obtained from primal
recovery:
\begin{equation}
 g_{\rm FW}(\Gamma)=\Tr[G(\Gamma-\widehat\Gamma)]
 =\lambda_P(\Gamma)-\lambda_D(\Gamma).
 \label{eq:cg-recovery-gap}
\end{equation}

Thus the Frank--Wolfe gap provides both the stopping criterion
for the reduced QP and the width of the recovered SOS interval.
Under the assumptions of Proposition~3.5,
\[
0 \le \phi(\Gamma)-\min\phi \le g_{\rm FW}(\Gamma).
\]

Before declaring numerical convergence, the implementation recomputes
the operator product at the final iterate, independently of its cached
linear combination, evaluates the global gap twice, and checks spectral
feasibility. These final evaluations use the same numerical backend but fresh operator products, so they directly test drift in the cache. Section~\ref{sec:protocol} describes the numerical settings.
Section~\ref{sec:recovery-experiment} checks the recovered SOS--SDP
objectives and Gram matrices against the original FCIDUMP integrals.

\section{Numerical experiments}
\label{sec:numerics}
The experiments address three questions: whether the reduced solutions
recover the original SDP objectives, how operator storage and corrective
updates affect computation, and what orbital dimensions the full-space
method can reach. We distinguish the molecular and FCIDUMP experiments from the
hydrogen-ring experiments.

\subsection{Inputs, accuracy, and timing}
\label{sec:protocol}
The GPU campaigns use eight NVIDIA H100 80 GB devices and FP64 arithmetic.
Geometry-based inputs use PySCF~\cite{sun2020pyscf},
GPU4PySCF~\cite{li2024gpu4pyscf}, and full symmetric orthogonalization.
FCIDUMP inputs supply the one- and two-electron coefficients directly.
Each DRSM solve starts from $\Gamma_0=I_N/2$ and uses at most 32
cached feasible points. The inner weight tolerance is $10^{-11}$ with
at most 2000 inner steps. The outer acceptance threshold is a fresh
absolute gap of $10^{-8}$ Ha and spectral-box violation at most $10^{-12}$.

\begin{table}[tbp]
\centering\small
\caption{Benchmark instances and sources of the input integrals or geometries.
$N$ is the number of spatial orbitals.}
\label{tab:campaigns}
\begin{tabularx}{\linewidth}{@{}lrr>{\raggedright\arraybackslash}X@{}}
\toprule
Campaign & Instances & $N$ & Data source \\
\midrule
FCIDUMP & 65 & 2--250 &
Electronic-structure integral collection associated with Low et
al.~\cite{low2025,rubin2025integrals} and catalytic-system integrals
from von Burg et al.~\cite{vonburg2021,vonburg2021integrals}. \\
Molecular set & 11 & 378--1404 &
Pentacene (Psi4)~\cite{psi4pentacene}, coronene (Nottoli et
al.)~\cite{nottoli2023geometries}, fullerene cages (FullereneNet and
Nanoten)~\cite{fullerenenet2025data,nanotenfullerenes}, and S12L/S30L
complexes~\cite{risthaus2013s12l,sure2015s30l,oterodelarozaRefdata}. \\
$\mathrm{C}_{720}$ cold start & 1 & 10,080 &
FullereneNet~\cite{fullerenenet2025data}. \\
Hydrogen rings & 19 & 4--10,000 &
Regular polygons generated in this work with nearest-neighbor spacing
1.25~\AA{} in STO-3G. The benchmark follows
Rubin et al.~\cite{rubin2026}. Integrals were generated with
PySCF/GPU4PySCF~\cite{sun2020pyscf,li2024gpu4pyscf}. \\
\bottomrule
\end{tabularx}
\end{table}

Table~\ref{tab:campaigns} counts experiments rather than distinct molecules
across campaigns. The FCIDUMP collection represents 65 distinct QPs.
The direct-integral screening parameter is $10^{-20}$ for the cached
molecular/C720 campaigns and $10^{-18}$ for the hydrogen-ring experiments.
FISTA measurements use their recorded configurations,
including a $10^{-5}$ target for the geometry sets.

Two quantities are reported separately: the feasible-point objective
$U=\phi(\Gamma)$ and the gap $g$, with lower endpoint $L=U-g$.
The Hamiltonian scalar is added only for total energies. For the hydrogen-ring runs, final evaluation uses the original integral operator after
numerical spectral-box projection. The stable gap formula is
\begin{equation}
 g=\sum_{i,j}(q_i^Tu_j)^2
 \bigl[(\mu_i)_+p_j+(-\mu_i)_+(1-p_j)\bigr],
 \label{eq:mf-complementarity-gap}
\end{equation}
where $\Gamma=U_\Gamma\operatorname{diag}(p_j)U_\Gamma^T$ and
$G=U_G\operatorname{diag}(\mu_i)U_G^T$, with columns $u_j$ and $q_i$, respectively. It equals the oracle gap in exact
arithmetic and avoids subtracting two large traces. The DRSM
molecular and FCIDUMP campaigns use $\Tr[G(\Gamma-\widehat\Gamma)]$.
All acceptance gaps use new operator evaluations rather than the cached
linear combination.

Timing columns retain their recorded boundaries. The DRSM
molecular/FCIDUMP solve clocks include final fresh checks. 
Input loading, matrix construction, and backend
initialization are kept separate from optimization. Table and figure captions
identify the convention used for each comparison.
The comparisons report implementation-level performance under
the configurations used in each campaign. FISTA uses gap targets
of $10^{-5}$ for the molecular set and $10^{-3}$ for C720,
whereas DRSM uses $10^{-8}$. Screening thresholds and operator backends follow the settings specified above
and in the corresponding experiment descriptions.

\subsection{FCIDUMP Hamiltonians}
\label{sec:cached-fcidump}

We evaluate DRSM and matrix-free FISTA on 65 distinct FCIDUMP quadratic models with
$2\leq N\leq250$ spatial orbitals. Appendix~\ref{app:fcidump instance}
lists the instance names. The source collections include the
electronic-structure integrals released with the work of Low et
al.~\cite{low2025,rubin2025integrals} and the computational-catalysis
integrals of von Burg et al.~\cite{vonburg2021,vonburg2021integrals},
which supply the I--XVIII defCD/highCD series.

DRSM uses eight GPUs, double precision, and the cold start
$\Gamma_0=I_N/2$. We assemble $H_2$ explicitly, distribute the rows
of $8H_2$ across the devices, and cache at most 32 feasible points.
The matrix-free FISTA uses the same FCIDUMP coefficients
and evaluates the quadratic operator by integral contraction.

All DRSM runs pass the numerical stopping tests. Both independently
recomputed final Frank--Wolfe gaps lie in $[0,10^{-8}]$~Ha, and all
eigenvalues of the final $\Gamma$ lie in $[-10^{-12},1+10^{-12}]$.
The runs require 7--207 outer updates and 0.1533--0.6872~s of solve
time. The DRSM solve clock includes the initial operator evaluation
and the final gap checks. It excludes FCIDUMP parsing, operator
assembly, preprocessing checks, and device loading. 
The FISTA iteration clock includes in-loop gap evaluations
and checkpoint writes. It excludes setup, initialization, and the
final gap check.

Table~\ref{tab:cached-fcidump-selected} compares nine representative
instances using the zero-core convention $E_{\mathrm{ecore}}=0$.
Both algorithms minimize
\[
  \phi_0(\Gamma)=2\svec(H_1)^T\gamma+4\gamma^T H_2\gamma.
\]
For DRSM we report $\phi_0$ at the saved final iterate. For matrix-free
FISTA we report the numerical lower endpoint,
denoted by $L_{\mathrm F}=\phi_0(\Gamma_{\mathrm F})-g_{\mathrm F}$,
where $\Gamma_{\mathrm F}$ and $g_{\mathrm F}$ are FISTA's final
iterate and Frank--Wolfe gap. Appendix~\ref{app:fcidump-results}
gives the complete results for all 65 distinct models in
Tables~\ref{tab:mf-fcidump} and~\ref{tab:cached-fcidump-all}.

\begin{table}[tbp]
\centering\small
\setlength{\tabcolsep}{3pt}
\caption{Selected FCIDUMP results on eight GPUs and reference energies
from Table~II of Low et al.~\cite{low2025}. Objectives, lower endpoints,
and gaps are in Hartree. Times are in seconds.}
\label{tab:cached-fcidump-selected}
\resizebox{\linewidth}{!}{%
\begin{tabular}{@{}lrrrrrrrr@{}}
\toprule
& & \multicolumn{3}{c}{DRSM}
  & \multicolumn{3}{c}{Matrix-free FISTA}
  & \multicolumn{1}{c}{Low et al.} \\
\cmidrule(lr){3-5}\cmidrule(lr){6-8}\cmidrule(l){9-9}
Instance & $N$ & $\phi_0$ & $g_{\mathrm D}$ & $t_{\mathrm D}$
               & $L_{\mathrm F}$ & $g_{\mathrm F}$ & $t_{\mathrm F}$
               & $E_{\mathrm{ref}}^{(0)}$ \\
\midrule
Fe2S2 & 20
  & -117.582263172 & 8.11e-09 & 0.2545
  & -117.582263183 & 6.19e-09 & 1.0057 & -117.582 \\
Fe4S4 & 36
  & -329.079184286 & 8.55e-09 & 0.5107
  & -329.079184298 & 6.06e-09 & 2.4708 & -329.079 \\
FeMoco54 & 54
  & -272.525404810 & 6.96e-09 & 0.4077
  & -272.525404822 & 8.93e-09 & 4.1004 & -272.524 \\
XVIII56 & 56
  & -296.944463511 & 3.94e-09 & 0.2628
  & -296.944463521 & 8.75e-09 & 2.3376 & -296.941 \\
P450X58 & 58
  & -424.552688575 & 7.08e-09 & 0.3155
  & -424.552688584 & 5.20e-09 & 2.9340 & -424.550 \\
FeMoco76 & 76
  & -1123.164456436 & 9.42e-09 & 0.6872
  & -1123.164456447 & 6.61e-09 & 3.6761 & -1123.300 \\
XVIII100$^{\dagger}$ & 100
  & -635.039589819 & 6.93e-09 & 0.3455
  & -635.039589829 & 6.36e-09 & 2.2898 & -635.134$^{*}$ \\
XVIII150$^{\dagger}$ & 150
  & -1355.454406010 & 9.19e-09 & 0.3690
  & -1355.454406023 & 7.49e-09 & 3.1770 & -1355.571$^{*}$ \\
XVIII250$^{\dagger}$ & 250
  & -3683.935458162 & 8.61e-09 & 0.5929
  & -3683.935458177 & 8.34e-09 & 10.4287 & --- \\
\bottomrule
\end{tabular}%
}
\par\smallskip
\begin{minipage}{\linewidth}\footnotesize
$g_{\mathrm D}$ is the final DRSM gap.
$t_{\mathrm D}$ is DRSM solve time, and $t_{\mathrm F}$ is FISTA
iteration time.
A dagger marks an input without a complete recorded numerical PSD check.
An asterisk marks the DFTHC-based estimates reported by Low et al.\
for XVIII100 and XVIII150 after their SDP solver failed to converge
within the allotted time. A dash indicates that Table~II of
Ref.~\cite{low2025} has no entry for XVIII250.
\end{minipage}
\end{table}

\subsection{Recovery in the original SDP coordinates}
\label{sec:recovery-experiment}

For each of nine instances, we recover the five Gram matrices of the
original SOS--SDP from the final QP iterate $\widehat\Gamma$ to verify feasibility. Using the
original FCIDUMP integrals, we evaluate the SOS--SDP objective and check
the coefficient equations and positive-semidefinite constraints.
Table~\ref{tab:original-sdp-recovery} compares the QP and SOS--SDP
objectives with the energies reported by Low et al.~\cite{low2025}.
All energies use $E_{\mathrm{ecore}}=0$.

The QP column gives $\phi_0(\widehat\Gamma)$. For the original SDP,
the constant coefficient equation gives the minimization objective
$\widehat\lambda=4\widehat q+
\operatorname{Tr}(Q_\uparrow)+\operatorname{Tr}(Q_\downarrow)$,
where $\widehat q=e^{\mathsf T}G'''e$ is evaluated from the lower-right
block of $G_{\mathrm{SF}}$. The SOS--SDP column gives this objective
in the electronic-energy convention,
$E_{\mathrm{ecore}}-\widehat\lambda=-\widehat\lambda$.
For $H_2\succeq0$ and $0\preceq\widehat\Gamma\preceq I$,
Proposition~\ref{prop:anytime-recovery} gives
\begin{equation}
  -\widehat\lambda \le E_{\mathrm{SOS}}^\star
  \le \phi_0(\widehat\Gamma),
  \qquad
  \phi_0(\widehat\Gamma)+\widehat\lambda
  = g_{\mathrm{FW}}(\widehat\Gamma).
  \label{eq:original-sdp-recovery-gap}
\end{equation}

Appendix~\ref{sec:lift} gives the corresponding dual multipliers and
slack matrices in the original SDP coordinates.
The recovery gives $D_\uparrow=D_\downarrow$ and $Q_\uparrow=Q_\downarrow$,
so one minimum eigenvalue suffices for each spin pair.
We report $r_{\mathrm{eq}}$ as the largest absolute
residual over the two-body, one-body (both spins), and constant
coefficient equations.

\begin{table*}[tbp]
\centering
\caption{QP and SOS--SDP objectives and feasibility checks in the
original SDP coordinates. The three energy columns use the zero-core
convention. The spectral columns give signed minimum eigenvalues.}
\label{tab:original-sdp-recovery}
\begingroup
\footnotesize
\setlength{\tabcolsep}{2pt}
\newcommand{\rsci}[2]{\ensuremath{#1\mathord{\times}10^{#2}}}
\begin{tabular*}{\textwidth}{@{\extracolsep{\fill}}lrrrrrrr@{}}
\toprule
& \multicolumn{3}{c}{Electronic energy (Ha)}
& \multicolumn{3}{c}{$\lambda_{\min}$, $\sigma\in\{\uparrow,\downarrow\}$}
& \\
\cmidrule(lr){2-4}\cmidrule(lr){5-7}
Instance & Ref.~\cite{low2025} & QP & SOS--SDP ($-\widehat\lambda$)
& $D_\sigma$ & $Q_\sigma$ & $G_{\mathrm{SF}}$ & $r_{\mathrm{eq}}$ \\
\midrule
Fe2S2
 & $-117.582$ & $-117.582263172$ & $-117.582263180$
 & \rsci{-8.61}{-18} & \rsci{-3.40}{-17}
 & \rsci{-3.24}{-16} & \rsci{4.26}{-14} \\
Fe4S4
 & $-329.079$ & $-329.079184286$ & $-329.079184295$
 & \rsci{-2.21}{-17} & \rsci{-6.15}{-17}
 & \rsci{-4.17}{-16} & \rsci{5.68}{-14} \\
FeMoco54
 & $-272.524$ & $-272.525404810$ & $-272.525404817$
 & \rsci{-6.45}{-17} & \rsci{-7.98}{-17}
 & \rsci{-6.35}{-16} & \rsci{1.14}{-13} \\
XVIII56
 & $-296.941$ & $-296.944463511$ & $-296.944463515$
 & \rsci{-1.69}{-16} & \rsci{-5.18}{-17}
 & \rsci{-1.03}{-15} & \rsci{1.78}{-15} \\
P450X58
 & $-424.550$ & $-424.552688575$ & $-424.552688583$
 & \rsci{-7.68}{-16} & \rsci{-4.53}{-17}
 & \rsci{-3.99}{-16} & \rsci{1.14}{-13} \\
FeMoco76
 & $-1123.300$ & $-1123.164456436$ & $-1123.164456446$
 & \rsci{-1.08}{-16} & \rsci{-1.52}{-16}
 & \rsci{-9.11}{-16} & \rsci{2.27}{-13} \\
XVIII100
 & $-635.134^{*}$ & $-635.039589819$ & $-635.039589825$
 & \rsci{-1.37}{-16} & \rsci{-7.36}{-17}
 & \rsci{-5.87}{-16} & \rsci{2.27}{-13} \\
XVIII150
 & $-1355.571^{*}$ & $-1355.454406010$ & $-1355.454406020$
 & \rsci{-2.96}{-16} & \rsci{-1.35}{-16}
 & \rsci{-9.81}{-13} & \rsci{2.27}{-13} \\
XVIII250
 & --- & $-3683.935458162$ & $-3683.935458171$
 & \rsci{-1.89}{-16} & \rsci{-9.72}{-17}
 & \rsci{-2.21}{-12} & \rsci{9.09}{-13} \\
\bottomrule
\end{tabular*}
\par\smallskip
\begin{minipage}{\textwidth}\raggedright
$^{*}$For XVIII100 and XVIII150, Low et al.\ used DFTHC estimates after
their SDP calculation failed to converge within the available runtime.
The dash indicates that XVIII250 does not appear in their Table~II.
\end{minipage}
\endgroup
\end{table*}

For FeMoco76, the SOS--SDP objective on the electronic-energy scale is $-1123.164456446$~Ha,
approximately $0.136$~Ha above the published value. Its largest equality
residual is $2.27\times10^{-13}$, and the smallest eigenvalue among the
five Gram matrices is $-9.11\times10^{-16}$. Fe2S2 and Fe4S4 agree with
the published energies to the three decimal places reported. For FeMoco54,
XVIII56, and P450X58, the SOS--SDP values lie below the published values by
approximately $0.0014$, $0.0035$, and $0.0027$~Ha, respectively.

\subsection{FISTA and DRSM on molecular geometries}
\label{sec:cached-results}
We compare matrix-free FISTA and DRSM on the same 11 molecular
geometries in the cc-pVDZ basis. Both methods use eight NVIDIA H100 80~GB GPUs.
FISTA targets an absolute gap of $10^{-5}$~Ha, and DRSM targets
$10^{-8}$~Ha.
Each DRSM solve starts from $\Gamma_0=I_N/2$.

The geometries comprise pentacene, coronene, fullerene cages, and S12L/S30L complexes in cc-pVDZ, using the source coordinates~\cite{psi4pentacene,nottoli2023geometries,fullerenenet2025data,nanotenfullerenes,risthaus2013s12l,sure2015s30l,oterodelarozaRefdata}. FISTA uses direct-integral matrix-free products for all 11 molecular
geometries. DRSM uses precomputed symmetric $H_2$ blocks for pentacene,
coronene, $\mathrm{C}_{36}$, $\mathrm{C}_{40}$, and $\mathrm{C}_{60}$.
These blocks are distributed across the eight GPUs. The recorded total
times include loading and validation of the blocks, whose construction
is performed beforehand. For the other six molecular geometries, DRSM uses direct-integral matrix-free products
with a screening threshold of $10^{-20}$. In both cases, the iterative solver
uses the same cached correction scheme. The eight GPUs perform the
operator work; the spectral oracle and feasible-point cache remain
primarily on GPU 0.

\begin{table}[tbp]
\centering
\small
\setlength{\tabcolsep}{4pt}
\caption{Matrix-free FISTA on 11 cc-pVDZ molecular geometries using eight GPUs
and an absolute gap target of $10^{-5}$ Ha. $E_{\mathrm{SOS}}$ is the
numerical lower endpoint including nuclear repulsion. $T_{\rm{iter}}$ is the iteration time.}
\label{tab:mf-large}
\begin{tabular}{@{}lrrrrr@{}}
\toprule
Instance & $N$ & $E_{\mathrm{SOS}}$ (Ha) & Gap (Ha) & Iterations & $T_{\rm{iter}}$(s) \\
\midrule
Pentacene & 378 & -888.103 & $6.20\times10^{-6}$ & 380 & 167.944 \\
Coronene & 396 & -965.799 & $6.10\times10^{-6}$ & 380 & 168.615 \\
$\mathrm{C}_{36}$ & 504 & -1429.662 & $4.29\times10^{-6}$ & 660 & 205.879 \\
$\mathrm{C}_{40}$ & 560 & -1588.496 & $9.76\times10^{-6}$ & 1820 & 641.646 \\
$\mathrm{C}_{60}$ & 840 & -2383.320 & $9.23\times10^{-6}$ & 2300 & 1854.408 \\
S12L 2\_COMPLEX1 & 982 & -2406.835 & $8.94\times10^{-6}$ & 1420 & 1265.944 \\
S30L 06\_di & 1018 & -2629.180 & $9.97\times10^{-6}$ & 2980 & 2640.849 \\
$\mathrm{C}_{76}$ & 1064 & -3021.982 & $8.04\times10^{-6}$ & 2420 & 3910.819 \\
$\mathrm{C}_{84}$ & 1176 & -3340.272 & $3.97\times10^{-6}$ & 940 & 2044.834 \\
$\mathrm{C}_{100}$ & 1400 & -3974.501 & $6.04\times10^{-6}$ & 1880 & 5291.802 \\
S12L 5\_COMPLEX1 & 1404 & -3452.803 & $8.06\times10^{-6}$ & 2020 & 2523.200 \\
\bottomrule
\end{tabular}
\end{table}

\begin{table}[tbp]
\centering\footnotesize
\setlength{\tabcolsep}{2.5pt}
\caption{Cold-start DRSM results. E denotes
explicit, GPU-resident symmetric $H_2$ blocks; MF denotes direct-integral
matrix-free evaluation. $k$ counts outer updates; each run uses $k+3$ Hessian products,
including initialization and two final checks. $E_{\rm iter}$ is the terminal
logged QP objective plus nuclear repulsion, rounded to six decimals;
the gap is independently re-evaluated. $T_{\rm solve}$
includes the initial product, iteration/checkpoint work and final checks;
$T_{\rm total}$ is the total elapsed time, including setup and
loading. Neither includes construction of previously stored $H_2$.}
\label{tab:cached-cold}
\begin{tabular}{@{}lrr r r r r r@{}}
\toprule
Instance & $N$ & Backend & $k$ & $E_{\rm iter}$ (Ha) & Gap (Ha) & $T_{\rm solve}$ (s) & $T_{\rm total}$ (s)\\
\midrule
Pentacene & 378 & E & 43 & -888.103040 & $9.80\times10^{-9}$ & 1.11 & 46.79 \\
Coronene & 396 & E & 35 & -965.798515 & $7.59\times10^{-9}$ & 1.07 & 49.92 \\
$\mathrm C_{36}$ & 504 & E & 59 & -1429.662344 & $9.54\times10^{-9}$ & 2.44 & 93.81 \\
$\mathrm C_{40}$ & 560 & E & 101 & -1588.496234 & $8.47\times10^{-9}$ & 4.80 & 143.81 \\
$\mathrm C_{60}$ & 840 & E & 100 & -2383.320263 & $6.74\times10^{-9}$ & 22.73 & 605.62 \\
S12L 2\_COMPLEX1 & 982 & MF & 76 & -2406.834759 & $6.94\times10^{-9}$ & 47.59 & 56.59 \\
S30L 06\_di & 1018 & MF & 99 & -2629.180224 & $7.08\times10^{-9}$ & 64.85 & 74.22 \\
$\mathrm C_{76}$ & 1064 & MF & 88 & -3021.981622 & $8.99\times10^{-9}$ & 111.31 & 120.81 \\
$\mathrm C_{84}$ & 1176 & MF & 75 & -3340.272315 & $9.04\times10^{-9}$ & 120.94 & 130.69 \\
$\mathrm C_{100}$ & 1400 & MF & 115 & -3974.500567 & $9.26\times10^{-9}$ & 255.34 & 266.43 \\
S12L 5\_COMPLEX1 & 1404 & MF & 118 & -3452.802598 & $8.64\times10^{-9}$ & 110.96 & 120.52 \\
\bottomrule
\end{tabular}
\end{table}

All 11 DRSM runs satisfy the $10^{-8}$ gap target in two final evaluations
using newly computed operator products. The spectral-box feasibility
checks also pass at tolerance $10^{-12}$. The largest observed violation
across these runs is approximately $5.8\times10^{-15}$; the number of
retained points at logged outer iterates never exceeds 10.
Across the 11 molecular geometries to both tables, the FISTA
lower endpoint and the DRSM final-iterate energy agree to three
decimal places in Hartree in every case. Summed over
these 11 instances, FISTA records 20,715.94~s of iteration time and
DRSM records 743.14~s of solve time, a ratio of 27.9 to one.

\subsection{A 10,080-orbital cold-start calculation on C720}
\label{sec:cached-c720}
\begin{figure}[tbp]
\centering
\includegraphics[width=\textwidth]{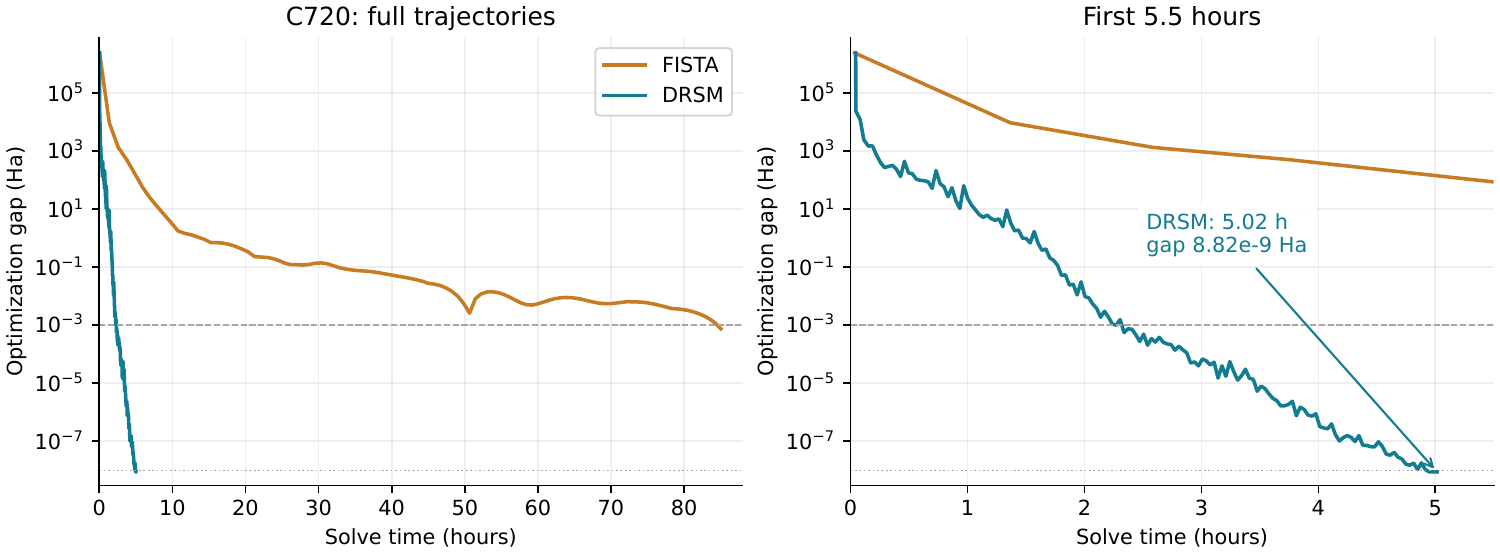}
\caption{C720 cold-start DRSM versus the eight-GPU
FISTA run. The right panel enlarges the first 5.5 hours of the same data.
DRSM reaches a freshly recomputed numerical gap of $8.82\times10^{-9}$ in 5.02 h;
FISTA terminates at about 85.04 h with gap
$7.55\times10^{-4}$. Their targets are $10^{-8}$ and $10^{-3}$,
respectively. Recorded clocks exclude model/backend setup; the cached
clock includes the initial action, checkpoints, and two final checks.}
\label{fig:cached-c720}
\end{figure}

We further solve the full spectral-box QP for $\mathrm C_{720}$ in the
cc-pVDZ basis, with $N=10{,}080$ and $d=50{,}808{,}240$. A dense FP64
$d$-by-$d$ matrix would require approximately 20.65 PB (decimal units).
We therefore use the direct-integral matrix-free backend with eight GPUs.
We obtained the unoptimized geometry of $\mathrm{C}_{720}$ from the
FullereneNet website~\cite{fullerenenet2025data}.
Figure~\ref{fig:cached-c720} shows the run initialized at $\Gamma_0=I_N/2$.

The run requires 147 outer updates and 150 Hessian products, including
the initial product and two fresh final evaluations. Its recorded solve
time is 18056.11 s (5.0156 h). The two final gaps are
$8.8181316\times10^{-9}$ and $8.8155536\times10^{-9}$ Ha. The final
matrix has minimum eigenvalue approximately $-5.36\times10^{-16}$ and
maximum eigenvalue approximately $1+6.0\times10^{-15}$. The final cache contains nine feasible points, with the $\mathsf V$ and $\mathsf W$ matrices occupying approximately
6.81~GiB. The calculation retains the full orbital basis and uses conventional four-center integral actions.

DRSM reaches a gap of $8.82\times10^{-9}$~Ha in 5.02~h,
compared with $7.55\times10^{-4}$~Ha in 85.04~h for the
FISTA. The convergence trajectories in Figure~\ref{fig:cached-c720}
show the improvement in both achieved accuracy and recorded
solve time.

\subsection{Full-space and cyclic-reduced hydrogen rings}
\label{sec:ring-experiment}
We compare full-space DRSM, cyclic-reduced DRSM, and BPSDP on regular hydrogen rings with special data pattern in the STO-3G basis with a nearest-neighbor spacing of 1.25~\AA. Both DRSM formulations start from \(\Gamma_0=I_N/2\) and target an absolute gap of \(10^{-8}\)~Ha. Full-space DRSM optimizes over the full matrix variable without imposing cyclic symmetry. Cyclic-reduced DRSM uses the symmetry reduction in Appendix~\ref{app:hydrogen-ring-reduction}. For each of the 19 ring sizes, we run full-space DRSM once and cyclic-reduced DRSM three times. The three reduced solves use the same precomputed reduced Hessian, and we report the median of their solve times.

For comparison, BPSDP in \texttt{libSDP} solves the bordered SDP~\eqref{eq:ordered-bordered} on six CPU threads. We accept a BPSDP run when its native primal residual, dual residual, and primal--dual objective gap meet their stopping criteria and an independent evaluation gives an absolute QP gap of at most \(10^{-5}\)~Ha. Appendix~\ref{app:full-hydrogen} reports the timing conventions and the per-instance results for all three methods.

\begin{figure}[tbp]
\centering
\includegraphics[width=.92\linewidth]{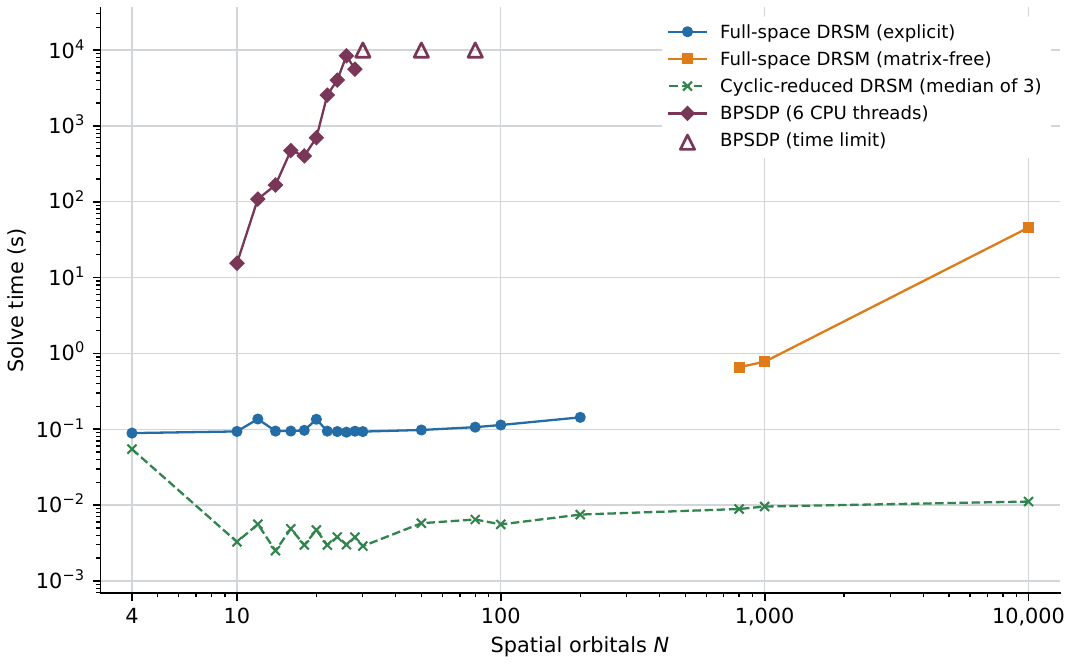}
\caption{Solve times for full-space DRSM, cyclic-reduced DRSM, and BPSDP on hydrogen rings. Each full-space DRSM point represents one solve. Each cyclic-reduced DRSM point is the median of three solve times for the same ring, excluding the one-time construction of the reduced Hessian. Open triangles mark BPSDP runs that reached the \(10{,}000\)~s process limit before meeting the acceptance criteria.}
\label{fig:hydrogen-new}
\end{figure}

Figure~\ref{fig:hydrogen-new} shows that both DRSM formulations reach the \(10^{-8}\)~Ha gap target for all 19 rings, up to \(N=10{,}000\). BPSDP meets its acceptance criteria for the ten tested sizes from \(N=10\) through \(N=28\) and reaches the process limit at \(N=30\), \(50\), and \(80\). Once the reduced Hessian has been constructed, cyclic symmetry lowers the cost of the operator and spectral calculations, with the largest solve-time savings at the largest ring sizes. Appendix~\ref{app:full-hydrogen} reports the complete per-instance results
for full-space DRSM, cyclic-reduced DRSM, and BPSDP.

\section{Discussion}
\label{sec:discussion}

The reformulation changes the computational structure of the SOS
problem as well as its dimension. In the expanded SDP, the generating
operators are represented through Gram matrices coupled by coefficient
matching. In the spectral-box QP, the feasible set has an explicit
linear minimization oracle, and the quadratic objective can be evaluated
through Coulomb contractions. These two properties connect the
optimization directly to spectral linear algebra and existing integral
engines. Feasible-point recovery preserves the connection to the
original certificate throughout the computation. This suggests a useful
principle for structured SOS optimization: eliminating redundant
coordinates is most valuable when it also exposes a feasible set and
operator representation that support inexpensive, specialized operations.

The cost of a Coulomb build also changes how optimization effort should
be allocated. When obtaining a new gradient is expensive, additional
work on a small corrective problem can be worthwhile if it makes
better use of operator products already available. DRSM implements
this allocation by retaining feasible points and their Hessian actions,
then optimizing their combination before requesting another action.
Its design therefore emphasizes the value obtained from each integral
build. The resulting trade-off is between the search space provided
by retained points and the memory and arithmetic needed to maintain
and combine them. An adaptive cache policy could use progress in the
restricted problem and the cost of a new action to decide when to
retain, remove, or generate directions. As integral calls become less
frequent, spectral decompositions and manipulation of the cached
matrices become increasingly important targets for parallelization.

The hydrogen rings illustrate how symmetry interacts with this
computational structure. For an invariant objective and a symmetric
initial point, the full-space algorithm can generate iterates that
retain the symmetry; the similar iteration counts observed in the
full and cyclic representations are consistent with this behavior.
Explicit symmetry reduction goes further by expressing those iterates
in coordinates that make their operator actions and spectral operations
cheaper. Thus symmetry can influence both the trajectory explored by
an algorithm and the cost of representing that trajectory. For less
symmetric molecular systems, an important question is which available
symmetries can be incorporated while preserving efficient Coulomb
contractions and a tractable spectral oracle. This offers a route to
combine chemical structure with the general optimization reduction.

Efficient solution of a fixed SOS model also makes it possible to
study the choice of model more systematically. The optimization gap
quantifies how accurately a chosen relaxation has been solved, while
the separation from a reference ground-state energy measures the
physical tightness of its bound. With the numerical optimization error
controlled, comparisons between generator spaces and symmetry
constraints can focus on the quality of the resulting bounds relative
to their computational cost. Stronger symmetry-adapted and higher-degree
SOS formulations are therefore natural subjects for further study.
The central structural question is whether they admit useful variable
eliminations, operator-based evaluations, and recovery maps analogous
to those developed here. Comparisons with FCI or DMRG reference
energies would then connect these computational developments to their
value for electronic-structure applications.

\section{Conclusion}

This work shows that the selected spin-free sum-of-squares energy
bound can be computed through an equivalent convex quadratic problem
over a spectral box, without directly solving the expanded
linear-objective SDP. The reformulation replaces the five PSD Gram matrices with a single
smaller symmetric matrix, reducing the variable count from quartic
to quadratic growth in the number of spatial orbitals while
preserving the optimal SOS bound. For $\mathrm{C}_{720}$ with
10,080 spatial orbitals, this reduces approximately
$2.06\times10^{16}$ scalar Gram-matrix variables to
$50{,}808{,}240$ independent scalar variables.
Explicit recovery maps connect every feasible QP iterate to primal
and dual witnesses in the original SOS formulation. The Frank--Wolfe
gap equals the separation between their objectives, providing both
a stopping criterion and an interval containing the optimal SOS energy.

The reformulated model also enables the development of a customized
solution method. DRSM combines a spectral linear minimization oracle
with cached Hessian actions and corrective optimization over a small
simplex of convex-combination weights. Reusing these actions allows
the inner weight updates to proceed without additional integral
evaluations. Together with an operator interface supporting both
stored coefficient blocks and GPU-accelerated matrix-free Coulomb
contractions, this design exploits the quadratic structure while
reducing repeated operator evaluations. On eight H100 GPUs, DRSM
solves the cold-start $\mathrm{C}_{720}$ instance with 10,080 spatial
orbitals to a freshly recomputed Frank--Wolfe gap of
$8.82\times10^{-9}$~Ha in 5.02 hours, demonstrating how the
reformulation and the customized solver together enable
large-scale SOS energy-bound calculations.

The hydrogen-ring experiments further demonstrate the computational
value of exploiting problem-specific data structure. Full-space DRSM
reaches 10,000 orbitals in 45.64 seconds of solve time, while explicit
cyclic reduction further decreases the optimization dimension and
the cost of operator and spectral calculations.
These results illustrate how exact reformulation, customized
optimization, and exploitation of additional structure can work
together to make large-scale SOS energy-bound calculations practical.

\section*{Data availability}
Numerical results are reported in the main text and appendices.
The source integral collections are identified by their version-specific
DOIs~\cite{rubin2025integrals,vonburg2021integrals}; geometry sources are
specified in Section~\ref{sec:cached-results} and
Section~\ref{sec:cached-c720}.
\bibliographystyle{plain}
\bibliography{references}
\clearpage
\appendix
\section{Hamiltonian reordering and exact coefficient matching}
\subsection{Reordering the electronic Hamiltonian}
\label{app:derivation of H}

The conventional state-dependent electronic Hamiltonian is
\[
 H=\sum_{ij,\sigma}t_{ij}a_{i\sigma}^{\dagger}a_{j\sigma}
 +\frac12\sum_{ijkl,\sigma,\tau}
 g_{ijkl}a_{i\sigma}^{\dagger}a_{k\tau}^{\dagger}a_{l\tau}a_{j\sigma}.
\]
The canonical anticommutation relations give
\begin{align*}
 E_{ij}E_{kl}
 &=\sum_{\sigma,\tau}a_{i\sigma}^{\dagger}a_{j\sigma}
                    a_{k\tau}^{\dagger}a_{l\tau}\\
 &=\delta_{jk}E_{il}
   -\sum_{\sigma,\tau}a_{i\sigma}^{\dagger}a_{k\tau}^{\dagger}
                     a_{j\sigma}a_{l\tau}\\
 &=\delta_{jk}E_{il}
   +\sum_{\sigma,\tau}a_{i\sigma}^{\dagger}a_{k\tau}^{\dagger}
                     a_{l\tau}a_{j\sigma}.
\end{align*}
Substitution yields
\[
 H=\sum_{il}\left(t_{il}-\frac12\sum_jg_{ijjl}\right)E_{il}
   +\frac12\sum_{ijkl}g_{ijkl}E_{ij}E_{kl}.
\]
Renaming the free indices gives \eqref{eq:h1h2}.
No constant is generated by this reordering.

\subsection{Proof of Proposition~\ref{prop:matching}}
\label{app:matching}
\begin{proof}
Let \(X_{ij}^\sigma=a_{i\sigma}^\dagger a_{j\sigma}\),
\(\beta_\sigma=\vecop(X^\sigma)\), and
\(\beta=\beta_\uparrow+\beta_\downarrow=\vecop(E)\).
Since \(x=P\beta\) and \(x^\dagger=\beta^T\), the expansion of the
five-block SOS has the form
\[
 \widehat H_{\rm SOS}
 =\beta^TAP\beta+
   \sum_{ij,\sigma}
     \bigl(D_\sigma-Q_\sigma+2(S+S^T)\bigr)_{ij}X_{ij}^\sigma
   +\kappa I,
 \qquad
 \kappa=4q+\Tr(Q_\uparrow)+\Tr(Q_\downarrow).
\]
The Hamiltonian is
\[
 H=\beta^TH_{\rm full}\beta+
       \sum_{ij,\sigma}(H_1)_{ij}X_{ij}^\sigma.
\]
Here \(A\succeq0\) because it is a principal block of
\(B^TG_{\rm SF}B\succeq0\).

To prove necessity, let \(\lvert\Omega\rangle\) be the vacuum,
\(\lvert i\uparrow\rangle=a_{i\uparrow}^\dagger\lvert\Omega\rangle\),
and
\(\lvert i\uparrow;k\downarrow\rangle
 =a_{i\uparrow}^\dagger a_{k\downarrow}^\dagger
   \lvert\Omega\rangle\), with the same creation order throughout.
For an operator \(\mathcal O\), define
\begin{align*}
 \mathcal C_{ij,kl}(\mathcal O)
 ={}&\langle i\uparrow;k\downarrow\vert\mathcal O
                         \vert j\uparrow;l\downarrow\rangle\\
 &-\delta_{kl}\langle i\uparrow\vert\mathcal O\vert j\uparrow\rangle
  -\delta_{ij}\langle k\downarrow\vert\mathcal O\vert l\downarrow\rangle
  +\delta_{ij}\delta_{kl}\langle\Omega\vert\mathcal O\vert\Omega\rangle .
\end{align*}
This functional annihilates constants and the even operators involving
only one spin sector. Anticommutation also gives
\[
 \mathcal C_{ij,kl}
       (X_{ab}^\uparrow X_{cd}^\downarrow)
       =\delta_{ia}\delta_{jb}\delta_{kc}\delta_{ld}.
\]
Indeed, on the one-up, one-down sector the product maps
\(\lvert j\uparrow;l\downarrow\rangle\) to
\(\delta_{bj}\delta_{dl}\lvert a\uparrow;c\downarrow\rangle\);
it vanishes on the vacuum and either one-electron sector.
The different-spin bilinears commute, so applying \(\mathcal C\) to
the claimed operator identity yields
\begin{equation}
 AP+PA=2H_{\rm full}.
 \label{eq:app-cross-matching}
\end{equation}
For any \(v\) with \(Pv=-v\), exchange symmetry gives
\(H_{\rm full}v=0\), and \eqref{eq:app-cross-matching} implies
\(v^TAv=0\).
Since \(A\succeq0\), this implies \(Av=0\).
Thus \(A\) annihilates the whole negative eigenspace of \(P\),
so \(AP=PA=A\). Equation~\eqref{eq:app-cross-matching} now gives
\(A=H_{\rm full}\).

The quadratic terms consequently cancel. Taking the vacuum
expectation of the remaining identity yields \(\lambda=\kappa\).
Taking matrix elements between one-electron states of each spin
then gives
\(D_\sigma-Q_\sigma+2(S+S^T)=H_1\) entry by entry.
This proves necessity.
Conversely, substituting \eqref{eq:exact-matching} into the expansion
immediately gives \(\widehat H_{\rm SOS}=H+\lambda I\).
The use of several particle-number sectors explains the full-Fock-space
assumption in the proposition.
\end{proof}

\subsection{Proof of Theorem~\ref{thm:pm-equivalence}}
\label{app:pm}
\begin{proof}
Given any original feasible point, congruence gives
$M=B^TG_{\rm SF}B\succeq0$.
Proposition~\ref{prop:matching} fixes the upper-left block to
$H_{\rm full}$ and gives exactly the remaining constraints and
objective of \eqref{eq:ordered-bordered}.
Keeping the four spin blocks therefore gives a feasible point with
the same objective.

For the reverse direction, $B[v;\alpha]=0$ implies $v=0$ from
its first block and then $\alpha e=0$. Since $N\geq1$,
$e\ne0$, so $\alpha=0$.
Thus $B$ has full column rank and $L=(B^TB)^{-1}B^T$ satisfies
$LB=I_{N^2+1}$.
Given a feasible bordered matrix $M$, set $G_{\rm SF}=L^TML$.
This is PSD, and
\[
 B^TG_{\rm SF}B=(LB)^TM(LB)=M.
\]
The maps $A,s,q$ are recovered exactly. Keeping the spin blocks and
setting $\lambda=4q+\Tr(Q_\uparrow)+\Tr(Q_\downarrow)$ now
recovers the original operator identity by
Proposition~\ref{prop:matching}, with the same objective.
The transformations therefore preserve the set of feasible
objective values and map optimal solutions to optimal solutions.
\end{proof}

\section{Proofs of compact reduction, duality, and recovery}
\label{app:sections-4-5}

\subsection{Proof of Theorem~\ref{thm:exact-reduction}}
\label{app:reduction}

\begin{proof}
We establish the equivalence by constructing
objective-preserving maps in both directions.
Recall that
\[
H_{\mathrm{full}}=UH_2U^T,
\qquad
U^TU=I_d.
\]
Define
\[
R=
\begin{bmatrix}
U&0\\
0&1
\end{bmatrix},
\qquad
R^TR=I_{d+1}.
\]

Given a feasible point of \eqref{eq:original-sos},
Theorem~\ref{thm:pm-equivalence} produces a feasible point
of \eqref{eq:ordered-bordered} with the same objective,
using
\[
M=B^TG_{\mathrm{SF}}B
=
\begin{bmatrix}
H_{\mathrm{full}}&s\\
s^T&q
\end{bmatrix}\succeq0.
\]
For any $v\in\ker H_{\mathrm{full}}$, the quadratic form
of $M$ on $[v;\alpha]$ is
\[
2\alpha v^Ts+\alpha^2q.
\]
Its nonnegativity for every $\alpha\in\mathbb R$
forces $v^Ts=0$.
Since $H_{\mathrm{full}}$ is symmetric,
\[
s\in(\ker H_{\mathrm{full}})^\perp
=\operatorname{range}(H_{\mathrm{full}})
\subseteq\operatorname{range}(U).
\]
Thus $s=Ut$ for a unique vector $t$, namely
$t=U^Ts$.
Congruence with $R^T$ now gives
\[
K:=R^TMR
=
\begin{bmatrix}
H_2&t\\
t^T&q
\end{bmatrix}\succeq0.
\]
Substituting $H_{\mathrm{full}}=UH_2U^T$
and $s=Ut$ also yields $M=RKR^T$.

Since $\operatorname{mat}(Ut)=\smat(t)$ is symmetric,
the one-body equations become
\[
D_\sigma-Q_\sigma
=H_1-4\smat(t),
\qquad \sigma\in\{\uparrow,\downarrow\}.
\]
Define
\[
D=\frac{D_\uparrow+D_\downarrow}{2},
\qquad
Q=\frac{Q_\uparrow+Q_\downarrow}{2}.
\]
Convexity of the positive-semidefinite cone gives
$D,Q\succeq0$, while averaging the two one-body equations
yields
\[
D-Q
=\frac12\bigl[
(D_\uparrow-Q_\uparrow)
+(D_\downarrow-Q_\downarrow)
\bigr]
=H_1-4\smat(t).
\]
By linearity of the trace,
\[
4q+\operatorname{Tr}(Q_\uparrow)
   +\operatorname{Tr}(Q_\downarrow)
=4q+2\operatorname{Tr}(Q).
\]
Thus $(t,q,D,Q)$ is feasible for \eqref{eq:p0}
and has the same objective value.

Conversely, let $(t,q,D,Q)$ be feasible for \eqref{eq:p0}.
Set
\[
D_\uparrow=D_\downarrow=D,
\qquad
Q_\uparrow=Q_\downarrow=Q,
\qquad
s=Ut.
\]
The duplicated matrices are positive semidefinite and satisfy
\[
D_\sigma-Q_\sigma
=D-Q=H_1-4\smat(t),
\qquad \sigma\in\{\uparrow,\downarrow\}.
\]
Their objective contribution is unchanged because
\[
\operatorname{Tr}(Q_\uparrow)
+\operatorname{Tr}(Q_\downarrow)
=2\operatorname{Tr}(Q).
\]
The bordered block
\[
K=
\begin{bmatrix}
H_2&t\\
t^T&q
\end{bmatrix}
\]
is positive semidefinite by feasibility of \eqref{eq:p0}.
Its congruence with $R$ therefore gives
\[
M=RKR^T
=
\begin{bmatrix}
H_{\mathrm{full}}&s\\
s^T&q
\end{bmatrix}\succeq0.
\]
Since $\operatorname{mat}(s)=\smat(t)$, the duplicated
spin variables also satisfy the one-body equations of
\eqref{eq:ordered-bordered}.
Hence this construction gives a feasible point of
\eqref{eq:ordered-bordered} with the same objective value.

The reverse map of Theorem~\ref{thm:pm-equivalence}
then recovers the original Gram variable as
\[
G_{\mathrm{SF}}
=L^TML
=L^TR
\begin{bmatrix}
H_2&t\\
t^T&q
\end{bmatrix}
R^TL,
\]
where $B$ and $L$ are defined in \eqref{eq:Bmap} and
\eqref{eq:left-inverse}.
In particular, the left-inverse relation $LB=I$ gives
\[
B^TG_{\mathrm{SF}}B
=(LB)^TM(LB)=M.
\]
The recovered Gram matrix is positive semidefinite,
and Theorem~\ref{thm:pm-equivalence} ensures that the
original coefficient equations and objective value
are preserved.
Therefore, all three formulations have the same set
of feasible objective values.
\end{proof}

\subsection{Proof of Theorem~\ref{thm:box-duality}}
\label{app:duality}

\begin{proof}
Let $H_2=L_2L_2^T$ be the factorization in
\eqref{eq:spectral-factorization}, where
$L_2\in\mathbb R^{d\times r}$ has full column rank and
$r=\operatorname{rank}(H_2)$.

We first verify the equivalence between \eqref{eq:p0}
and the square-root formulation
\eqref{eq:spectral-factor-primal}.
Suppose
\[
\begin{bmatrix}
H_2&t\\
t^T&q
\end{bmatrix}\succeq0.
\]
For any $v\in\ker H_2$, nonnegativity of the associated
quadratic form on $(v,a)$ for all $a\in\mathbb R$
forces $v^Tt=0$.
Hence
\[
t\in\operatorname{range}(H_2)
=\operatorname{range}(L_2),
\]
so $t=L_2y$ for a unique $y\in\mathbb R^r$.
Writing
\[
W=
\begin{bmatrix}
L_2&0\\
0&1
\end{bmatrix},
\]
we have
\[
\begin{bmatrix}
H_2&t\\
t^T&q
\end{bmatrix}
=
W
\begin{bmatrix}
I_r&y\\
y^T&q
\end{bmatrix}
W^T.
\]
Since $W$ has full column rank, congruence with a left
inverse of $W$ shows that the bordered matrix is PSD
if and only if the middle matrix is PSD.
The ordinary Schur complement therefore gives
\[
q\geq\|y\|_2^2.
\]
This proves \eqref{eq:generalized-schur} without requiring
$H_2$ to be invertible.

Every feasible point of \eqref{eq:p0} consequently gives
a feasible point $(y,D,Q)$ of
\eqref{eq:spectral-factor-primal} with objective
\[
4\|y\|_2^2+2\Tr(Q)
\leq 4q+2\Tr(Q).
\]
Conversely, every feasible $(y,D,Q)$ of
\eqref{eq:spectral-factor-primal} gives a feasible point
of \eqref{eq:p0} through
\[
t=L_2y,\qquad q=\|y\|_2^2,
\]
with exactly the same objective value.
Thus the two formulations have the same infimum.

The square-root formulation is strictly feasible.
Indeed, choose
\[
\alpha>\max\{0,-\lambda_{\min}(H_1)\},
\qquad
y=0,\qquad
Q=\alpha I_N,\qquad
D=H_1+\alpha I_N.
\]
Then $D,Q\succ0$, and
\[
D-Q=H_1-4\smat(L_2y).
\]
Its objective is finite, and every feasible objective
value is nonnegative.

We next compute the Lagrange dual.
Use the symmetric equality multiplier $2\Gamma$ and
write $\gamma=\svec(\Gamma)$.
The Lagrangian in \eqref{eq:spectral-lagrangian} is
\[
\begin{aligned}
\mathcal L(y,D,Q;\Gamma)
={}&4\|y\|_2^2
+8\gamma^TL_2y
-2\Tr(H_1\Gamma)\\
&+2\Tr(\Gamma D)
+2\Tr\bigl((I_N-\Gamma)Q\bigr).
\end{aligned}
\]
Its infimum over $D,Q\succeq0$ is finite exactly when
\[
0\preceq\Gamma\preceq I_N.
\]
To see necessity, if $\Gamma$ has a negative eigenvector
$v$, choosing $D=\rho vv^T$ and letting
$\rho\to+\infty$ makes $\Tr(\Gamma D)\to-\infty$.
The same argument applies to $Q$ when $I_N-\Gamma$
has a negative eigenvalue.
Conversely, if both coefficient matrices are PSD,
the two trace terms are nonnegative and attain zero
at $D=Q=0$.

For a spectral-box feasible $\Gamma$, completing the
square gives
\[
\begin{aligned}
4\|y\|_2^2+8\gamma^TL_2y
&=4\|y+L_2^T\gamma\|_2^2
  -4\gamma^TH_2\gamma.
\end{aligned}
\]
Thus the infimum over $y$ is attained at
\[
y_\Gamma=-L_2^T\gamma,
\]
and the dual function equals
\[
-4\gamma^TH_2\gamma
-2\Tr(H_1\Gamma)
=-\phi(\Gamma).
\]
The Lagrange dual is therefore
\[
\max_{0\preceq\Gamma\preceq I_N}
\{-\phi(\Gamma)\}.
\]

The convex Slater condition applies to the strictly
feasible square-root formulation and yields zero
duality gap.
Since the spectral box is compact and $\phi$ is
continuous, the dual optimum is attained. Consequently,
\[
\lambda_*
=\max_{0\preceq\Gamma\preceq I_N}\{-\phi(\Gamma)\}
=-\min_{0\preceq\Gamma\preceq I_N}\phi(\Gamma),
\]
which proves \eqref{eq:boxdual} and attainment for
\eqref{eq:phi}.

It remains to establish primal attainment.
Choose $\beta$ equal to the objective value of the
strictly feasible point constructed above.
On the feasible sublevel set
\[
4\|y\|_2^2+2\Tr(Q)\leq\beta,
\]
we have
\[
\|y\|_2^2\leq\frac{\beta}{4},
\qquad
\|Q\|_F\leq\Tr(Q)\leq\frac{\beta}{2},
\]
where the second inequality uses $Q\succeq0$.
The equality constraint
\[
D=Q+H_1-4\smat(L_2y)
\]
then also bounds $D$. More explicitly, the isometry
of $\smat$ gives
\[
\|D\|_F
\leq\|Q\|_F+\|H_1\|_F
+4\|L_2\|_2\|y\|_2.
\]
This sublevel set is nonempty, closed, and bounded,
and hence compact.
The continuous objective therefore attains its
minimum on it, proving attainment for
\eqref{eq:spectral-factor-primal}.
Finally, the map
\[
t=L_2y,\qquad q=\|y\|_2^2
\]
recovers an optimizer of \eqref{eq:p0}.
When $r=0$, $H_2=0$ and $\phi(\Gamma)=2\Tr(H_1\Gamma)$.
The primal choice $t=0$, $q=0$, $D=(H_1)_+$, $Q=(-H_1)_+$
attains $2\Tr((-H_1)_+)$; the spectral-box minimum is its negative.
This establishes the same duality and attainment statements in this case.
\end{proof}

\subsection{Proof of Proposition~\ref{prop:anytime-recovery}}
\label{app:proof-anytime-recovery}

\begin{proof}
Fix $0\preceq\Gamma\preceq I_N$ and write $\gamma=\svec(\Gamma)$.
By \eqref{eq:feasible-recovery},
$t_\Gamma=-H_2\gamma$ and $q_\Gamma=\gamma^T H_2\gamma$.
Since $H_2=H_2^T\succeq0$, the bordered block in \eqref{eq:p0}
satisfies
\[
\begin{aligned}
 \begin{bmatrix}H_2&t_\Gamma\\ t_\Gamma^T&q_\Gamma\end{bmatrix}
 &=\begin{bmatrix}
     H_2&-H_2\gamma\\
     -\gamma^T H_2&\gamma^T H_2\gamma
   \end{bmatrix} \\
 &=\begin{bmatrix}I_d\\-\gamma^T\end{bmatrix}
   H_2\begin{bmatrix}I_d&-\gamma\end{bmatrix}
 \succeq0 .
\end{aligned}
\]
The recovered matrices $D_\Gamma=(K_\Gamma)_+$ and
$Q_\Gamma=(-K_\Gamma)_+$ are positive semidefinite.
Using the spectral identity $K=K_+-(-K)_+$ for a symmetric matrix $K$
gives the remaining equality in \eqref{eq:p0},
\[
 D_\Gamma-Q_\Gamma
 =K_\Gamma
 =H_1+4\smat(H_2\gamma)
 =H_1-4\smat(t_\Gamma).
\]
Thus $(t_\Gamma,q_\Gamma,D_\Gamma,Q_\Gamma)$ is feasible for
\eqref{eq:p0}, with objective $\lambda_P(\Gamma)=4q_\Gamma+2\Tr(Q_\Gamma)$.
The objective-preserving lift in Theorem~\ref{thm:exact-reduction}
yields an original SOS certificate with this same objective.
Since $0\preceq\Gamma\preceq I_N$, its dual objective is
$\lambda_D(\Gamma)=-\phi(\Gamma)$ by Theorem~\ref{thm:box-duality}.

Write $K_\Gamma=O\operatorname{diag}(\kappa_1,\ldots,\kappa_N)O^T$.
Then $G=2K_\Gamma$ has eigenvalues $2\kappa_i$ in the same basis.
For any $0\preceq Z\preceq I_N$, the diagonal entries
$z_i=(O^T ZO)_{ii}$ satisfy $0\le z_i\le1$. Hence
\[
 \Tr(GZ)=2\sum_{i=1}^N\kappa_i z_i
 \ge 2\sum_{i:\,\kappa_i<0}\kappa_i.
\]
The projector
$\widehat\Gamma=O\operatorname{diag}(\mathbf{1}_{\{\kappa_i<0\}})O^T$
belongs to the spectral box and attains this bound. Since the
eigenvalues of $Q_\Gamma=(-K_\Gamma)_+$ are $\max\{-\kappa_i,0\}$,
\[
 \min_{0\preceq Z\preceq I_N}\Tr(GZ)
 =\Tr(G\widehat\Gamma)
 =2\sum_{i:\,\kappa_i<0}\kappa_i
 =-2\Tr(Q_\Gamma).
\]

For the other trace in the gap, the isometry of $\svec$ gives
\[
 \Tr\!\left(\smat(H_2\gamma)\Gamma\right)
 =(H_2\gamma)^T\gamma
 =q_\Gamma.
\]
It follows from $G=2H_1+8\smat(H_2\gamma)$ and
$\phi(\Gamma)=2\Tr(H_1\Gamma)+4q_\Gamma$ that
\[
 \Tr(G\Gamma)=2\Tr(H_1\Gamma)+8q_\Gamma
              =\phi(\Gamma)+4q_\Gamma .
\]
Combining the two trace identities yields
\[
\begin{aligned}
 \Tr\!\bigl(G(\Gamma-\widehat\Gamma)\bigr)
 &=\Tr(G\Gamma)-\Tr(G\widehat\Gamma)\\
 &=\phi(\Gamma)+4q_\Gamma+2\Tr(Q_\Gamma)\\
 &=\lambda_P(\Gamma)-\lambda_D(\Gamma).
\end{aligned}
\]
The left-hand side equals the Frank--Wolfe gap because
$\widehat\Gamma$ minimizes $\Tr(GZ)$ over the spectral box.

By weak duality,
$\lambda_D(\Gamma)\le\lambda_*\le\lambda_P(\Gamma)$.
Subtracting from $E_{\rm ecore}$ reverses these inequalities and gives
\[
 E_{\rm ecore}-\lambda_P(\Gamma)
 \le E_{\rm SOS}^{\star}
 = E_{\rm ecore}-\lambda_*
 \le E_{\rm ecore}-\lambda_D(\Gamma).
\]
Finally, $\lambda_D(\Gamma)=-\phi(\Gamma)$ and
$\lambda_P(\Gamma)=\lambda_D(\Gamma)+g_{\rm FW}(\Gamma)$.
Substituting these expressions at the two endpoints yields
\[
 E_{\rm ecore}+\phi(\Gamma)-g_{\rm FW}(\Gamma)
 \le E_{\rm SOS}^{\star}
 \le E_{\rm ecore}+\phi(\Gamma),
\]
as claimed.
\end{proof}

\subsection{Proof of Proposition~\ref{prop:recovery}}
\label{app:recovery}

\begin{proof}
Let $\Gamma^\star$ solve \eqref{eq:phi}, and let
$(t^\star,q^\star,D^\star,Q^\star)$ be defined by
\eqref{eq:primal-recovery}. Applying
Proposition~\ref{prop:anytime-recovery} at $\Gamma=\Gamma^\star$
gives a feasible solution of \eqref{eq:p0} and an original SOS
certificate with the same objective value.

Set $G^\star=\nabla\phi(\Gamma^\star)$. By first-order optimality,
\[
 \Tr\!\bigl(G^\star(Z-\Gamma^\star)\bigr)\geq 0
 \qquad\text{for every }0\preceq Z\preceq I_N.
\]
Thus $\Gamma^\star$ minimizes $\Tr(G^\star Z)$ over the spectral box,
and hence
\[
 g_{\rm FW}(\Gamma^\star)
 =\Tr(G^\star\Gamma^\star)
  -\min_{0\preceq Z\preceq I_N}\Tr(G^\star Z)
 =0.
\]
The gap identity in \eqref{eq:energy-interval} now gives
$\lambda_P(\Gamma^\star)=\lambda_D(\Gamma^\star)$.
Since weak duality places $\lambda_*$ between these two values,
\[
 4q^\star+2\Tr(Q^\star)
 =\lambda_P(\Gamma^\star)
 =\lambda_*
 =\lambda_D(\Gamma^\star)
 =-\phi(\Gamma^\star).
\]
Therefore the recovered solution is optimal for \eqref{eq:p0},
and its objective-preserving lift is optimal for the original
SOS--SDP. The energy identities follow from
$\varepsilon_{\rm SOS}^{\star}=-\lambda_*=\phi(\Gamma^\star)$
and $E_{\rm SOS}^{\star}=E_{\rm ecore}+\phi(\Gamma^\star)$.
\end{proof}

\section{The dual witness in the original SDPA coordinates}
\label{sec:lift}
We express the spin-free constraints from Appendix~F of
Ref.~\cite{low2025} in the sparse SDPA format. We order the equality
constraints by placing the two one-body groups first and the two-body
group last. With this ordering, let $\xi=\operatorname{vec}(\Gamma)$. The ordering and
normalization below specify the files used in our recovery test.
The equality multiplier vector is
\begin{equation}
 y=\begin{bmatrix}-\xi\\-\xi\\-4(\xi\otimes\xi)\end{bmatrix}.
 \label{eq:original-dual-multipliers}
\end{equation}
Substitution in the original dual slack equations gives
\begin{equation}
\begin{aligned}
 Z_{D_\uparrow}=Z_{D_\downarrow}&=\Gamma,\\
 Z_{Q_\uparrow}=Z_{Q_\downarrow}&=I_N-\Gamma,\\
 Z_{G_{\rm SF}}&=4vv^T,\qquad
 v=\begin{bmatrix}\vecop(\Gamma)\\\vecop(I_N-\Gamma)\end{bmatrix}.
\end{aligned}
\label{eq:original-dual-slacks}
\end{equation}
The small slacks are PSD for every spectral-box feasible point; the
large slack is given directly as a PSD outer product. The original dual
objective is
\begin{equation}
 b^Ty=-2\Tr(H_1\Gamma)-4\gamma^TH_2\gamma=-\phi(\Gamma).
 \label{eq:original-dual-value}
\end{equation}
Together with the primal lift in Theorem~\ref{thm:exact-reduction}, these
formulas represent both witnesses using the orbital-space iterate.
The spectral-box solver works directly with $\Gamma$. The
$N^4+2N^2$ multipliers can be formed from $\Gamma$ when an explicit
dual witness in SDPA coordinates is required.
\section{An additional cyclic reduction for regular hydrogen rings}
\label{app:hydrogen-ring-reduction}
For a regular ring with one identical real $s$-type basis function per
atom, cyclic symmetry permits a further reduction of the spectral-box QP.
This is additional to the general SOS-to-QP reformulation.

\begin{proposition}[Cyclic spectral coordinates]
\label{prop:ring-reduction}
Assume positive-definite AO overlap, full symmetric orthogonalization,
and $H_2\succeq0$. For $m=\lfloor N/2\rfloor+1$, set
\[
 w_k=\begin{cases}1,&k=0\text{ or }(N\text{ even and }k=N/2),\\
 2,&\text{otherwise},\end{cases}\qquad 0\leq k<m,
\]
and, with orbital indices $0\leq i,j<N$, define
\begin{equation}
 (E_k)_{ij}=\frac{\sqrt{w_k}}{N}
 \cos\!\left(\frac{2\pi k(j-i)}{N}\right),\quad
 B_{\rm cyc}=[\svec(E_0),\ldots,\svec(E_{m-1})].
 \label{eq:ring-basis}
\end{equation}
Then $B_{\rm cyc}^TB_{\rm cyc}=I_m$, and the full QP has the same
optimal value as
\begin{equation}
 \min_{0\leq z_k\leq\sqrt{w_k}}
 \frac12 z^TA_{\rm cyc}z+b_{\rm cyc}^Tz,\qquad
 A_{\rm cyc}=8B_{\rm cyc}^TH_2B_{\rm cyc},\quad
 b_{\rm cyc}=2B_{\rm cyc}^T\svec(H_1).
 \label{eq:ring-reduced-qp}
\end{equation}
An optimizer recovers a full-space optimizer through
$\Gamma=\sum_k z_kE_k$.
\end{proposition}
\begin{proof}
Let $T$ cyclically permute the AO functions. Rotation of the ring
preserves the overlap and all integral coefficients under this permutation.
Since $S_{\rm AO}^{-1/2}$ commutes with $T$, the same action is valid in
the symmetrically orthogonalized orbital basis. Consequently,
$\phi(T\Gamma T^T)=\phi(\Gamma)$, and orthogonal conjugation preserves
the spectral box. The cyclic average
\[
 \overline\Gamma=\frac1N\sum_{j=0}^{N-1}T^j\Gamma(T^j)^T
\]
is feasible, invariant, and has objective no greater than that of
$\Gamma$, by convexity. Applying this to an optimizer gives an invariant
optimizer.

The invariant symmetric matrices are real symmetric circulant matrices.
Fourier diagonalization~\cite{Gray2006} identifies their independent
eigenvalues $p_k$, of multiplicities $w_k$. The corresponding orthogonal
projectors are $\sqrt{w_k}E_k$, so the $E_k$ are Frobenius-orthonormal.
Writing $z_k=\sqrt{w_k}p_k$ gives
$\Gamma=\sum_k z_kE_k$ and converts $0\preceq\Gamma\preceq I_N$ into
$0\leq z_k\leq\sqrt{w_k}$. Substitution into \eqref{eq:phi} gives
\eqref{eq:ring-reduced-qp}. Every reduced feasible point is full-space
feasible with the same objective, and an invariant optimizer exists.
\end{proof}
For $N=10{,}000$, this changes 50,005,000 packed coordinates to 5001
cyclic coordinates. The full-space experiment in
Section~\ref{sec:ring-experiment} does not impose this parametrization.

\clearpage
\section{FCIDUMP instance identifiers}
\label{app:fcidump instance}

Instances are ordered by parsed spatial-orbital count and then original
instance name. The identifiers below correspond to
Table~\ref{tab:mf-fcidump} and Table~\ref{tab:cached-fcidump-all}.

\begingroup
\fontsize{8.5}{9}\selectfont
\setlength{\tabcolsep}{4pt}
\begin{longtable}{@{}lrp{0.76\linewidth}@{}}
\caption{FCIDUMP instance identifiers. $N$ is the parsed number
of spatial orbitals.}\label{tab:mf-names}\\
\toprule
ID & $N$ & Original instance name \\
\midrule
\endfirsthead
\toprule
ID & $N$ & Original instance name \\
\midrule
\endhead
\bottomrule
\endfoot
F01 & 2 & \path|defCD__VIII-cas2-fb-2e2o| \\
F02 & 2 & \path|highCD__VIII-highCD-cas2-fb-2e2o| \\
F03 & 4 & \path|Hchain__h2-8so-full| \\
F04 & 4 & \path|defCD__VIII-IX-cas4-fb-4e4o| \\
F05 & 4 & \path|defCD__XVIII-cas4-fb-4e4o| \\
F06 & 4 & \path|highCD__VIII-IX-highCD-cas4-fb-4e4o| \\
F07 & 4 & \path|highCD__XVIII-highCD-cas4-fb-4e4o| \\
F08 & 5 & \path|defCD__I-cas5-fb-4e5o| \\
F09 & 5 & \path|highCD__I-highCD-cas5-fb-4e5o| \\
F10 & 6 & \path|defCD__II-III-cas6-fb-8e6o| \\
F11 & 6 & \path|defCD__II-cas6-fb-8e6o| \\
F12 & 6 & \path|highCD__II-III-highCD-cas6-fb-8e6o| \\
F13 & 6 & \path|highCD__II-highCD-cas6-fb-8e6o| \\
F14 & 8 & \path|Hchain__h4-16so-full| \\
F15 & 11 & \path|defCD__V-cas11-fb-12e11o| \\
F16 & 11 & \path|highCD__V-highCD-cas11-fb-12e11o| \\
F17 & 12 & \path|Hchain__h6-24so-full| \\
F18 & 16 & \path|Hchain__h8-32so-full| \\
F19 & 16 & \path|defCD__I-cas5-fb-14e16o| \\
F20 & 16 & \path|defCD__IX-cas16-fb-16e16o| \\
F21 & 16 & \path|highCD__I-highCD-cas5-fb-14e16o| \\
F22 & 16 & \path|highCD__IX-highCD-cas16-fb-16e16o| \\
F23 & 20 & \path|Fe2S2| \\
F24 & 20 & \path|defCD__XVIII-cas4-fb-28e20o| \\
F25 & 20 & \path|highCD__XVIII-highCD-cas4-fb-28e20o| \\
F26 & 23 & \path|defCD__VIII-IX-cas4-fb-36e23o| \\
F27 & 23 & \path|highCD__VIII-IX-highCD-cas4-fb-36e23o| \\
F28 & 24 & \path|defCD__V-cas11-fb-32e24o| \\
F29 & 24 & \path|highCD__V-highCD-cas11-fb-32e24o| \\
F30 & 26 & \path|defCD__II-cas6-fb-34e26o| \\
F31 & 26 & \path|defCD__IX-cas16-fb-32e26o| \\
F32 & 26 & \path|highCD__II-highCD-cas6-fb-34e26o| \\
F33 & 26 & \path|highCD__IX-highCD-cas16-fb-32e26o| \\
F34 & 29 & \path|defCD__II-III-cas6-fb-38e29o| \\
F35 & 29 & \path|defCD__VIII-cas2-fb-40e29o| \\
F36 & 29 & \path|highCD__II-III-highCD-cas6-fb-38e29o| \\
F37 & 29 & \path|highCD__VIII-highCD-cas2-fb-40e29o| \\
F38 & 33 & \path|defCD__II-hf-fb-38e33o| \\
F39 & 36 & \path|Fe4S4| \\
F40 & 52 & \path|defCD__I-cas5-fb-48e52o| \\
F41 & 52 & \path|highCD__I-highCD-cas5-fb-48e52o| \\
F42 & 54 & \path|FeMoco54| \\
F43 & 56 & \path|XVIII56| \\
F44 & 56 & \path|highCD__XVIII-highCD-cas4-fb-64e56o| \\
F45 & 58 & \path|P450X58| \\
F46 & 59 & \path|defCD__VIII-IX-cas4-fb-72e59o| \\
F47 & 59 & \path|highCD__VIII-IX-highCD-cas4-fb-72e59o| \\
F48 & 60 & \path|defCD__V-cas11-fb-68e60o| \\
F49 & 60 & \path|highCD__V-highCD-cas11-fb-68e60o| \\
F50 & 62 & \path|defCD__II-cas6-fb-70e62o| \\
F51 & 62 & \path|defCD__IX-cas16-fb-68e62o| \\
F52 & 62 & \path|highCD__II-highCD-cas6-fb-70e64o| \\
F53 & 62 & \path|highCD__IX-highCD-cas16-fb-68e62o| \\
F54 & 65 & \path|defCD__II-III-cas6-fb-74e65o| \\
F55 & 65 & \path|defCD__VIII-cas2-fb-76e65o| \\
F56 & 65 & \path|highCD__II-III-highCD-cas6-fb-74e65o| \\
F57 & 65 & \path|highCD__VIII-highCD-cas2-fb-76e65o| \\
F58 & 71 & \path|defCD__II-hf-fb-74e71o| \\
F59 & 76 & \path|FeMoco76| \\
F60 & 100 & \path|XVIII100| \\
F61 & 100 & \path|highCD__XVIII-highCD-cas4-fb-100e100o| \\
F62 & 150 & \path|XVIII150| \\
F63 & 150 & \path|highCD__XVIII-highCD-cas4-fb-150e150o| \\
F64 & 250 & \path|XVIII250| \\
F65 & 250 & \path|highCD__XVIII-highCD-cas4-fb-250e250o| \\
\end{longtable}
\endgroup

\section{FISTA and DRSM results for the FCIDUMP collection}
\label{app:fcidump-results}
\label{app:legacy-fista}\label{app:mf-experiments}\label{app:cached-fcidump}

Tables~\ref{tab:mf-fcidump} and~\ref{tab:cached-fcidump-all} report
results for the same 65 distinct FCIDUMP instances, identified in
Table~\ref{tab:mf-names}. Both methods use eight GPUs and an absolute
gap target of $10^{-8}$~Ha. FISTA applies the quadratic operator by
matrix-free contraction. DRSM uses an explicitly assembled $H_2$.

Both FISTA and DRSM minimize the same zero-core objective $\phi_0$.
Table~\ref{tab:mf-fcidump} reports the
numerical lower endpoint
$E_{\mathrm{SOS}}^{(0)}
=\phi_0(\Gamma_{\mathrm{F}})-g_{\mathrm{FW}}(\Gamma_{\mathrm{F}})$.
Table~\ref{tab:cached-fcidump-all} reports the objective at the final DRSM
iterate and the larger of two independently recomputed final gaps.
Both tables exclude the FCIDUMP scalar $E_{\mathrm{ecore}}$. FISTA times cover
the iteration loop, while DRSM solve times include the initial operator
evaluation and the final gap checks. Section~\ref{sec:protocol} gives
the timing conventions.

\begin{table}[htbp]
\centering
\fontsize{8}{10}\selectfont
\setlength{\tabcolsep}{3pt}
\caption{Matrix-free FISTA results for 65 distinct FCIDUMP instances
on eight GPUs, with a uniform absolute gap target of $10^{-8}$ Ha.
The energy column gives the lower endpoint with the FCIDUMP scalar removed.
Energies and gaps are in Hartree, and times are in seconds.}
\label{tab:mf-fcidump}
\begin{tabular}{@{}lrrrr@{\hspace{5mm}}lrrrr@{}}
\toprule
\multicolumn{5}{c}{Instances F01--F33}
& \multicolumn{5}{c}{Instances F34--F65} \\
\cmidrule(lr){1-5}\cmidrule(lr){6-10}
ID & $N$ & $E_{\mathrm{SOS}}^{(0)}$ & gap & t(s)
& ID & $N$ & $E_{\mathrm{SOS}}^{(0)}$ & gap & t(s) \\
\midrule
F01 & 2 & -2.075 & $9.78\times10^{-17}$ & 0.3067
& F34 & 29 & -168.430 & $1.05\times10^{-10}$ & 0.4747 \\
F02 & 2 & -2.075 & $2.10\times10^{-18}$ & 0.2806
& F35 & 29 & -178.528 & $1.20\times10^{-10}$ & 0.5071 \\
F03 & 4 & -1.862 & $2.22\times10^{-12}$ & 0.4316
& F36 & 29 & -168.426 & $1.58\times10^{-10}$ & 0.4773 \\
F04 & 4 & -5.090 & $2.01\times10^{-10}$ & 0.2794
& F37 & 29 & -178.526 & $1.20\times10^{-10}$ & 0.4847 \\
F05 & 4 & -6.733 & $1.22\times10^{-9}$ & 0.2994
& F38 & 33 & -160.143 & $2.43\times10^{-10}$ & 0.5464 \\
F06 & 4 & -5.090 & $2.01\times10^{-10}$ & 0.3184
& F39 & 36 & -329.079 & $6.06\times10^{-9}$ & 2.4708 \\
F07 & 4 & -6.733 & $1.22\times10^{-9}$ & 0.3024
& F40 & 52 & -165.382 & $5.26\times10^{-9}$ & 2.3053 \\
F08 & 5 & -5.090 & $9.16\times10^{-15}$ & 0.3294
& F41 & 52 & -165.383 & $8.60\times10^{-9}$ & 2.3146 \\
F09 & 5 & -5.090 & $9.64\times10^{-15}$ & 0.3716
& F42 & 54 & -272.525 & $8.93\times10^{-9}$ & 4.1004 \\
F10 & 6 & -17.113 & $4.59\times10^{-13}$ & 0.3526
& F43 & 56 & -296.944 & $8.75\times10^{-9}$ & 2.3376 \\
F11 & 6 & -18.264 & $6.42\times10^{-9}$ & 0.2852
& F44 & 56 & -296.941 & $5.71\times10^{-9}$ & 2.7533 \\
F12 & 6 & -17.113 & $4.59\times10^{-13}$ & 0.3541
& F45 & 58 & -424.553 & $5.20\times10^{-9}$ & 2.9340 \\
F13 & 6 & -18.264 & $6.41\times10^{-9}$ & 0.2858
& F46 & 59 & -387.187 & $6.28\times10^{-9}$ & 1.8932 \\
F14 & 8 & -4.977 & $1.65\times10^{-9}$ & 0.4118
& F47 & 59 & -387.188 & $9.98\times10^{-9}$ & 1.8244 \\
F15 & 11 & -39.210 & $4.71\times10^{-9}$ & 0.3761
& F48 & 60 & -357.723 & $6.54\times10^{-9}$ & 3.2432 \\
F16 & 11 & -39.209 & $4.71\times10^{-9}$ & 0.3401
& F49 & 60 & -357.717 & $9.56\times10^{-9}$ & 3.2733 \\
F17 & 12 & -8.630 & $1.39\times10^{-9}$ & 0.5053
& F50 & 62 & -363.920 & $9.32\times10^{-9}$ & 1.9201 \\
F18 & 16 & -11.974 & $2.38\times10^{-9}$ & 0.6863
& F51 & 62 & -369.214 & $4.87\times10^{-9}$ & 2.0231 \\
F19 & 16 & -31.859 & $8.80\times10^{-9}$ & 0.3534
& F52 & 62 & -363.894 & $8.71\times10^{-9}$ & 1.9634 \\
F20 & 16 & -59.424 & $6.46\times10^{-12}$ & 1.0406
& F53 & 62 & -369.214 & $5.76\times10^{-9}$ & 2.0149 \\
F21 & 16 & -31.859 & $8.79\times10^{-9}$ & 0.3510
& F54 & 65 & -404.752 & $3.16\times10^{-10}$ & 3.2292 \\
F22 & 16 & -59.424 & $6.45\times10^{-12}$ & 0.4334
& F55 & 65 & -420.514 & $4.83\times10^{-10}$ & 1.6541 \\
F23 & 20 & -117.582 & $6.19\times10^{-9}$ & 1.0057
& F56 & 65 & -404.748 & $1.64\times10^{-9}$ & 3.2541 \\
F24 & 20 & -91.306 & $3.47\times10^{-9}$ & 0.4315
& F57 & 65 & -420.511 & $5.40\times10^{-9}$ & 1.9615 \\
F25 & 20 & -91.303 & $3.47\times10^{-9}$ & 0.4118
& F58 & 71 & -400.259 & $8.35\times10^{-9}$ & 3.1923 \\
F26 & 23 & -152.726 & $5.10\times10^{-9}$ & 0.9108
& F59 & 76 & -1123.164 & $6.61\times10^{-9}$ & 3.6761 \\
F27 & 23 & -152.726 & $4.47\times10^{-9}$ & 0.9264
& F60 & 100 & -635.040 & $6.36\times10^{-9}$ & 2.2898 \\
F28 & 24 & -145.429 & $1.08\times10^{-9}$ & 0.4871
& F61 & 100 & -635.039 & $1.03\times10^{-9}$ & 2.2829 \\
F29 & 24 & -145.420 & $5.79\times10^{-10}$ & 0.4997
& F62 & 150 & -1355.454 & $7.49\times10^{-9}$ & 3.1770 \\
F30 & 26 & -140.951 & $3.68\times10^{-9}$ & 0.5636
& F63 & 150 & -1355.456 & $4.93\times10^{-9}$ & 3.7253 \\
F31 & 26 & -152.425 & $6.96\times10^{-10}$ & 0.5372
& F64 & 250 & -3683.935 & $8.34\times10^{-9}$ & 10.4287 \\
F32 & 26 & -140.924 & $6.00\times10^{-9}$ & 0.4264
& F65 & 250 & -3683.933 & $8.37\times10^{-10}$ & 8.9644 \\
F33 & 26 & -152.425 & $6.97\times10^{-10}$ & 0.5762
&  &  &  &  &  \\
\bottomrule
\end{tabular}
\par\smallskip
\begin{minipage}{\linewidth}
\footnotesize
Energies are rounded to three decimal places after removing the input scalar.
Gaps and stopping tests use unrounded values. Time includes the iteration
loop, in-loop gap evaluations and checkpoint writes. It excludes preparation,
loading, initialization and the final check. Full instance names appear in
Appendix~\ref{app:fcidump instance}.
\end{minipage}
\end{table}

\begin{table}[tbp]
\centering
\fontsize{8}{10}\selectfont
\setlength{\tabcolsep}{2.5pt}
\caption{DRSM results for the 65 distinct FCIDUMP quadratic programs.
The feasible-point objective $\phi_0$ excludes the input scalar.
Objectives and fresh gaps are in Hartree.}
\label{tab:cached-fcidump-all}
\resizebox{\linewidth}{!}{%
\begin{tabular}{@{}lrrrrr@{\hspace{4mm}}lrrrrr@{}}
\toprule
\multicolumn{6}{c}{Instances F01--F33}
& \multicolumn{6}{c}{Instances F34--F65} \\
\cmidrule(lr){1-6}\cmidrule(lr){7-12}
ID & $N$ & $\phi_0$ & Fresh gap & Updates & $t$ (s)
& ID & $N$ & $\phi_0$ & Fresh gap & Updates & $t$ (s) \\
\midrule
F01 & 2 & -2.074791316 & 3.70e-09 & 11 & 0.1546
& F34 & 29 & -168.430202277 & 7.48e-09 & 17 & 0.1792 \\
F02 & 2 & -2.074831224 & 3.68e-09 & 11 & 0.1565
& F35 & 29 & -178.528231143 & 3.48e-09 & 16 & 0.1823 \\
F03 & 4 & -1.861729768 & 7.66e-09 & 15 & 0.1780
& F36 & 29 & -168.425719333 & 7.40e-09 & 17 & 0.1877 \\
F04 & 4 & -5.090051779 & 2.30e-10 & 7 & 0.1533
& F37 & 29 & -178.526183343 & 3.52e-09 & 16 & 0.1759 \\
F05 & 4 & -6.733282605 & 5.10e-09 & 9 & 0.1544
& F38 & 33 & -160.142532946 & 5.51e-09 & 16 & 0.1864 \\
F06 & 4 & -5.090025919 & 2.30e-10 & 7 & 0.1538
& F39 & 36 & -329.079184286 & 8.55e-09 & 149 & 0.5107 \\
F07 & 4 & -6.733170648 & 5.11e-09 & 9 & 0.1607
& F40 & 52 & -165.382362324 & 8.44e-09 & 54 & 0.2519 \\
F08 & 5 & -5.089799126 & 4.66e-09 & 9 & 0.1633
& F41 & 52 & -165.382736023 & 6.82e-09 & 55 & 0.2737 \\
F09 & 5 & -5.089808373 & 4.74e-09 & 9 & 0.1557
& F42 & 54 & -272.525404810 & 6.96e-09 & 109 & 0.4077 \\
F10 & 6 & -17.113078625 & 3.97e-09 & 13 & 0.1743
& F43 & 56 & -296.944463511 & 3.94e-09 & 57 & 0.2628 \\
F11 & 6 & -18.263753904 & 6.53e-09 & 18 & 0.1655
& F44 & 56 & -296.940977959 & 4.11e-09 & 53 & 0.2491 \\
F12 & 6 & -17.113003309 & 3.97e-09 & 13 & 0.1701
& F45 & 58 & -424.552688575 & 7.08e-09 & 73 & 0.3155 \\
F13 & 6 & -18.263746959 & 6.75e-09 & 18 & 0.1771
& F46$^{\dagger}$ & 59 & -387.186688461 & 9.97e-09 & 36 & 0.2216 \\
F14 & 8 & -4.977471229 & 6.58e-09 & 12 & 0.1720
& F47$^{\dagger}$ & 59 & -387.187827935 & 5.50e-09 & 37 & 0.2176 \\
F15 & 11 & -39.210186753 & 2.33e-09 & 17 & 0.1889
& F48$^{\dagger}$ & 60 & -357.722550477 & 5.32e-09 & 48 & 0.2639 \\
F16 & 11 & -39.209234566 & 2.33e-09 & 17 & 0.1730
& F49$^{\dagger}$ & 60 & -357.717180053 & 7.58e-09 & 52 & 0.2512 \\
F17 & 12 & -8.630417480 & 7.22e-09 & 19 & 0.1895
& F50$^{\dagger}$ & 62 & -363.919526163 & 8.86e-09 & 62 & 0.2721 \\
F18 & 16 & -11.973618436 & 3.83e-09 & 31 & 0.2043
& F51$^{\dagger}$ & 62 & -369.214178305 & 4.67e-09 & 38 & 0.2375 \\
F19 & 16 & -31.859246496 & 8.79e-09 & 12 & 0.1813
& F52$^{\dagger}$ & 62 & -363.893609739 & 3.11e-09 & 47 & 0.2639 \\
F20 & 16 & -59.424305266 & 9.77e-09 & 10 & 0.1575
& F53$^{\dagger}$ & 62 & -369.213645417 & 9.91e-09 & 37 & 0.2209 \\
F21 & 16 & -31.859397651 & 3.25e-09 & 13 & 0.1633
& F54$^{\dagger}$ & 65 & -404.752481272 & 6.62e-09 & 47 & 0.2695 \\
F22 & 16 & -59.424056862 & 9.72e-09 & 10 & 0.1796
& F55$^{\dagger}$ & 65 & -420.513953469 & 8.60e-09 & 38 & 0.2424 \\
F23 & 20 & -117.582263172 & 8.11e-09 & 49 & 0.2545
& F56$^{\dagger}$ & 65 & -404.748145159 & 6.71e-09 & 54 & 0.3244 \\
F24 & 20 & -91.305847839 & 3.90e-09 & 14 & 0.1771
& F57$^{\dagger}$ & 65 & -420.511239806 & 7.91e-09 & 37 & 0.2309 \\
F25 & 20 & -91.303491506 & 3.91e-09 & 14 & 0.1826
& F58$^{\dagger}$ & 71 & -400.258675363 & 7.50e-09 & 57 & 0.3278 \\
F26 & 23 & -152.725567118 & 3.74e-09 & 32 & 0.2155
& F59 & 76 & -1123.164456436 & 9.42e-09 & 207 & 0.6872 \\
F27 & 23 & -152.726387291 & 3.62e-09 & 32 & 0.2133
& F60$^{\dagger}$ & 100 & -635.039589819 & 6.93e-09 & 71 & 0.3455 \\
F28 & 24 & -145.428912292 & 8.10e-09 & 18 & 0.1828
& F61$^{\dagger}$ & 100 & -635.039331474 & 9.17e-09 & 72 & 0.3572 \\
F29 & 24 & -145.419977914 & 8.05e-09 & 18 & 0.2316
& F62$^{\dagger}$ & 150 & -1355.454406010 & 9.19e-09 & 78 & 0.3690 \\
F30 & 26 & -140.950621568 & 7.16e-09 & 16 & 0.2381
& F63$^{\dagger}$ & 150 & -1355.456063359 & 7.89e-09 & 72 & 0.3544 \\
F31 & 26 & -152.424990250 & 3.16e-09 & 17 & 0.1792
& F64$^{\dagger}$ & 250 & -3683.935458162 & 8.61e-09 & 83 & 0.5929 \\
F32 & 26 & -140.924091153 & 6.95e-09 & 16 & 0.1754
& F65$^{\dagger}$ & 250 & -3683.932591172 & 8.10e-09 & 74 & 0.5249 \\
F33 & 26 & -152.425239932 & 3.15e-09 & 17 & 0.1805
&  &  &  &  &  &  \\
\bottomrule
\end{tabular}%
}
\par\smallskip
\begin{minipage}{\linewidth}\footnotesize
$t$ is solve-only time. A dagger marks an instance without a complete
recorded numerical check that $H_2\succeq0$.
\end{minipage}
\end{table}

\section{Complete hydrogen-ring results}
\label{app:full-hydrogen}

The three tables give the measurements plotted in
Figure~\ref{fig:hydrogen-new}. All rings have nearest-neighbor spacing
1.25~\AA{} and use the STO-3G basis.
The primal and dual columns give the SDP minimization objectives
$\lambda_P$ and $\lambda_D$ in Hartree, with the Hamiltonian scalar
excluded. For DRSM, $\lambda_D=-\phi(\Gamma)$ and
$\lambda_P=\lambda_D+g_{\mathrm{FW}}(\Gamma)$.

\subsection*{Full-space DRSM}

Each row is one cold-start solve on eight GPUs, starting from
$\Gamma_0=I_N/2$ in $d=N(N+1)/2$ dimensions. All 19 runs meet the
$10^{-8}$~Ha gap target. The gap is the larger of two original-integral
evaluations after spectral-box feasibility cleanup. Solve time excludes
operator construction, initialization, and the two final checks.

\begingroup
\small
\setlength{\tabcolsep}{3.5pt}
\begin{longtable}{@{}rrcrrrrr@{}}
\caption{Full-space DRSM. E denotes explicitly assembled $H_2$
and MF denotes direct-integral matrix-free evaluation.}
\label{tab:hydrogen-all}\\
\toprule
$N$ & $d$ & Backend & Updates & Solve time (s)
& Primal & Dual & Gap (Ha) \\
\midrule
\endfirsthead
\toprule
$N$ & $d$ & Backend & Updates & Solve time (s)
& Primal & Dual & Gap (Ha) \\
\midrule
\endhead
\bottomrule\endfoot
4 & 10 & E & 3 & 0.088992 & 4.541 & 4.541 & $4.15\times10^{-14}$ \\
10 & 55 & E & 1 & 0.093606 & 16.375 & 16.375 & $2.00\times10^{-15}$ \\
12 & 78 & E & 4 & 0.136375 & 20.548 & 20.548 & $2.28\times10^{-14}$ \\
14 & 105 & E & 1 & 0.094550 & 25.021 & 25.021 & $1.33\times10^{-15}$ \\
16 & 136 & E & 4 & 0.094681 & 29.478 & 29.478 & $2.15\times10^{-14}$ \\
18 & 171 & E & 2 & 0.096834 & 34.146 & 34.146 & $3.95\times10^{-15}$ \\
20 & 210 & E & 4 & 0.135326 & 38.819 & 38.819 & $4.52\times10^{-14}$ \\
22 & 253 & E & 2 & 0.094402 & 43.647 & 43.647 & $2.22\times10^{-14}$ \\
24 & 300 & E & 3 & 0.093257 & 48.490 & 48.490 & $7.93\times10^{-15}$ \\
26 & 351 & E & 2 & 0.091507 & 53.453 & 53.453 & $4.86\times10^{-14}$ \\
28 & 406 & E & 3 & 0.094247 & 58.439 & 58.439 & $1.18\times10^{-13}$ \\
30 & 465 & E & 2 & 0.093150 & 63.519 & 63.519 & $5.50\times10^{-15}$ \\
50 & 1,275 & E & 5 & 0.097826 & 116.769 & 116.769 & $1.31\times10^{-13}$ \\
80 & 3,240 & E & 6 & 0.106551 & 202.815 & 202.815 & $2.04\times10^{-13}$ \\
100 & 5,050 & E & 5 & 0.113767 & 262.984 & 262.984 & $1.43\times10^{-14}$ \\
200 & 20,100 & E & 7 & 0.143554 & 584.683 & 584.683 & $8.65\times10^{-14}$ \\
800 & 320,400 & MF & 9 & 0.658373 & 2808.218 & 2808.218 & $3.14\times10^{-12}$ \\
1000 & 500,500 & MF & 10 & 0.776495 & 3604.724 & 3604.724 & $7.00\times10^{-13}$ \\
10000 & 50,005,000 & MF & 11 & 45.639181 & 45793.747 & 45793.747 & $2.17\times10^{-10}$ \\
\end{longtable}
\endgroup

Across the 19 full-space runs, the total elapsed time including
initialization and validation was 283.47~s. GPU memory was sampled
every 2~s. The largest samples over the full queue were 28,695~MiB
on GPU~0 and 11,355~MiB on each of the other GPUs. The operator work
was distributed across the GPUs, while GPU~0 held the primary direction
cache and evaluated the spectral oracle.

\subsection*{Cyclic-reduced DRSM}

The cyclic reduction in Appendix~\ref{app:hydrogen-ring-reduction}
gives $m=\lfloor N/2\rfloor+1$ coordinates. For each size, we performed
three cold-start solves from the reduced coordinates of $\Gamma_0=I_N/2$.
All three produced the same final coordinates. The reported solve time
excludes construction of the reduced operator. For each size, we expanded
the final coordinates of the first run into the full matrix and evaluated
the gap twice using the original integral operator. The table reports the
larger of these two gaps.

\begingroup
\small
\setlength{\tabcolsep}{3.5pt}
\begin{longtable}{@{}rrrrrrr@{}}
\caption{Cyclic-reduced DRSM. Times are medians of three solves
with the same precomputed reduced Hessian.}
\label{tab:hydrogen-cyclic-all}\\
\toprule
$N$ & $m$ & Updates & Median solve time (s)
& Primal & Dual & Gap (Ha) \\
\midrule
\endfirsthead
\toprule
$N$ & $m$ & Updates & Median solve time (s)
& Primal & Dual & Gap (Ha) \\
\midrule
\endhead
\bottomrule\endfoot
4 & 3 & 3 & 0.054678 & 4.541 & 4.541 & $7.162\times10^{-16}$ \\
10 & 6 & 1 & 0.003278 & 16.375 & 16.375 & $1.280\times10^{-16}$ \\
12 & 7 & 4 & 0.005584 & 20.548 & 20.548 & $5.745\times10^{-15}$ \\
14 & 8 & 1 & 0.002505 & 25.021 & 25.021 & $7.966\times10^{-16}$ \\
16 & 9 & 4 & 0.004866 & 29.478 & 29.478 & $3.865\times10^{-15}$ \\
18 & 10 & 2 & 0.002978 & 34.146 & 34.146 & $1.632\times10^{-15}$ \\
20 & 11 & 4 & 0.004698 & 38.819 & 38.819 & $3.949\times10^{-15}$ \\
22 & 12 & 2 & 0.002970 & 43.647 & 43.647 & $1.627\times10^{-14}$ \\
24 & 13 & 3 & 0.003796 & 48.490 & 48.490 & $7.422\times10^{-15}$ \\
26 & 14 & 2 & 0.003004 & 53.453 & 53.453 & $4.829\times10^{-15}$ \\
28 & 15 & 3 & 0.003772 & 58.439 & 58.439 & $2.476\times10^{-14}$ \\
30 & 16 & 2 & 0.002879 & 63.519 & 63.519 & $5.430\times10^{-15}$ \\
50 & 26 & 5 & 0.005769 & 116.769 & 116.769 & $1.058\times10^{-14}$ \\
80 & 41 & 6 & 0.006460 & 202.815 & 202.815 & $5.924\times10^{-15}$ \\
100 & 51 & 5 & 0.005565 & 262.984 & 262.984 & $7.209\times10^{-15}$ \\
200 & 101 & 7 & 0.007498 & 584.683 & 584.683 & $3.442\times10^{-14}$ \\
800 & 401 & 9 & 0.008877 & 2808.218 & 2808.218 & $2.896\times10^{-13}$ \\
1000 & 501 & 10 & 0.009578 & 3604.724 & 3604.724 & $8.158\times10^{-13}$ \\
10000 & 5001 & 11 & 0.011134 & 45793.747 & 45793.747 & $4.125\times10^{-12}$ \\
\end{longtable}
\endgroup

For $N=10{,}000$, the $5001\times5001$ reduced Hessian occupied
190.81~MiB and took approximately 3.40~h to construct.

\subsection*{BPSDP}

The boundary-point solver in \texttt{libSDP} solves the bordered
SDP~\eqref{eq:ordered-bordered} on six CPU threads. Acceptance requires
the native primal residual, dual residual, and primal--dual objective gap
to pass their stopping tests, together with an independently evaluated
absolute QP gap of at most $10^{-5}$~Ha. Solve times exclude problem
construction, solver setup, and the final QP gap evaluation. The two
objective columns give the native BPSDP primal and dual values. The
last column gives the independently evaluated QP gap.

\begingroup
\small
\setlength{\tabcolsep}{3.5pt}
\begin{longtable}{@{}rlrrrr@{}}
\caption{BPSDP on six CPU threads.}
\label{tab:hydrogen-bpsdp-all}\\
\toprule
$N$ & Status & Solve time (s)
& Primal & Dual & QP gap (Ha) \\
\midrule
\endfirsthead
\toprule
$N$ & Status & Solve time (s)
& Primal & Dual & gap (Ha) \\
\midrule
\endhead
\bottomrule\endfoot
10 & Accepted & 15.45 & 16.375345550 & 16.375339057 & $2.42\times10^{-6}$ \\
12 & Accepted & 108.38 & 20.547870388 & 20.547880387 & $8.18\times10^{-6}$ \\
14 & Accepted & 166.15 & 25.020971152 & 25.020981149 & $4.99\times10^{-7}$ \\
16 & Accepted & 472.40 & 29.478434333 & 29.478435162 & $2.75\times10^{-6}$ \\
18 & Accepted & 401.62 & 34.146435216 & 34.146441510 & $7.18\times10^{-7}$ \\
20 & Accepted & 697.36 & 38.818582350 & 38.818591540 & $4.46\times10^{-6}$ \\
22 & Accepted & 2548.10 & 43.646703208 & 43.646710288 & $1.23\times10^{-6}$ \\
24 & Accepted & 4020.07 & 48.489804882 & 48.489814232 & $1.10\times10^{-6}$ \\
26 & Accepted & 8401.98 & 53.453408934 & 53.453409483 & $9.97\times10^{-6}$ \\
28 & Accepted & 5600.36 & 58.438601254 & 58.438591767 & $4.32\times10^{-7}$ \\
30 & Time limit & -- & -- & -- & -- \\
50 & Time limit & -- & -- & -- & -- \\
80 & Time limit  & -- & -- & -- & -- \\
\end{longtable}
\endgroup

Dashes indicate that the process reached the 10,000~s limit before
meeting the acceptance criteria. Figure~\ref{fig:hydrogen-new} places
all three open triangles at this limit.
\end{document}